%% file: main.tex
\documentclass[hidelinks,onefignum,onetabnum]{siamart220329}

\input{ex_shared}

\ifpdf
\hypersetup{
  pdftitle={},
  pdfauthor={}
}
\fi

\usepackage[colorinlistoftodos]{todonotes}
\usepackage{amsmath}
\usepackage{standalone}
\usepackage{tikz}
\usepackage{amsfonts} %
\usepackage[normalem]{ulem}
\usepackage{marginnote}
\usepackage{pgfplots}
\usepackage{placeins}
\usepackage{soul}
\usepackage{amssymb}

\usepackage{graphbox}

\newcommand{\mjb}[1]{\textcolor{red!70!red}{MJB: #1}}
\newcommand{\onehalf}{\frac{1}{2}}
\DeclareMathSymbol{\sm}{\mathbin}{AMSa}{"39}
\usetikzlibrary{patterns.meta,perspective}
\usepackage{caption} 
\usepackage{atveryend}
\usepackage{pdfpages}
\AtVeryEndDocument{\includepdf[pages={-}]{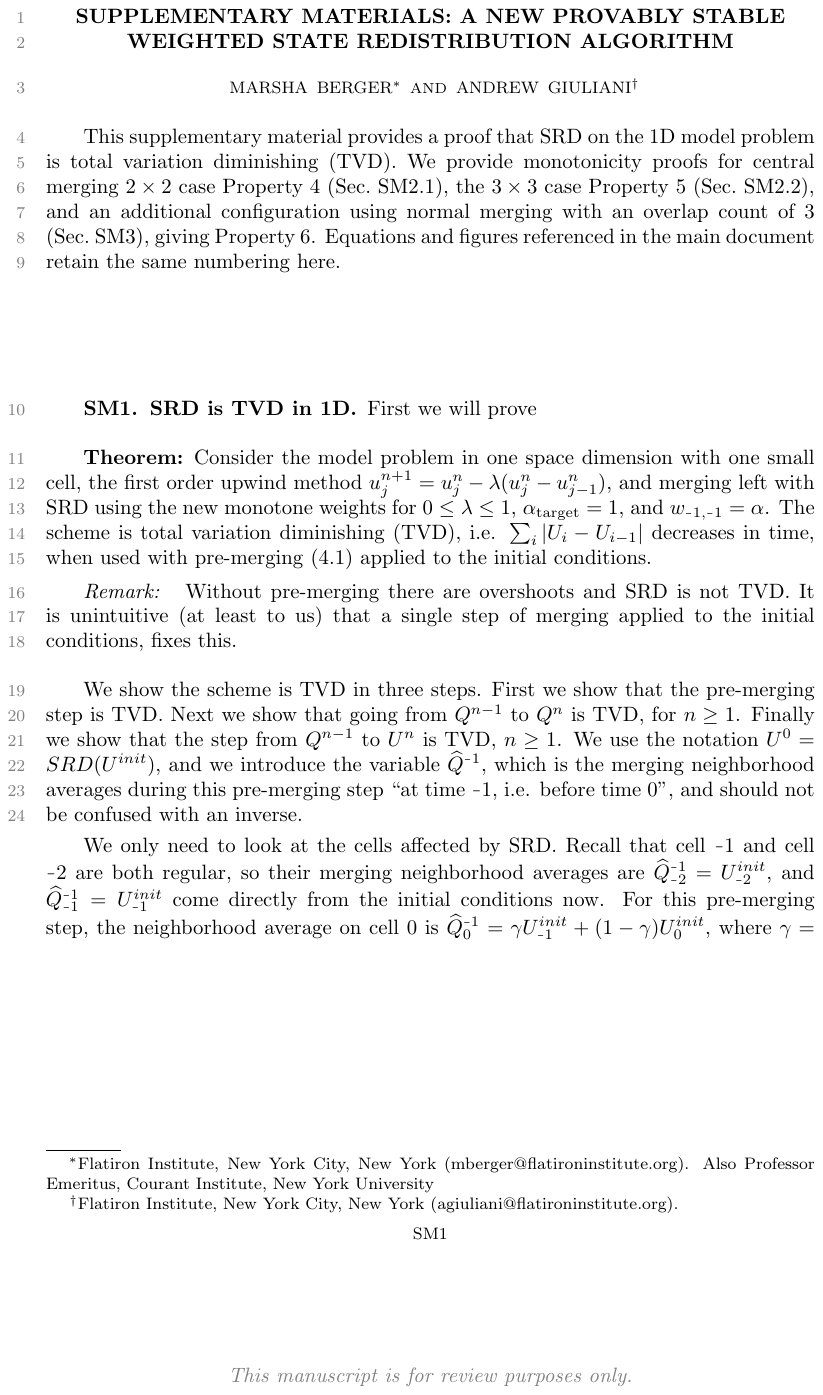}}

\begin{document}

\maketitle

\begin{abstract} 
We propose a practical finite volume method on cut cells using state redistribution.
Our algorithm is provably monotone, total variation diminishing, and GKS stable in many situations, and shuts off continuously as the cut cell size approaches a target value.  
Our analysis reveals why original state redistribution works so well: it results in a monotone scheme for most configurations, though at times subject to a slightly smaller CFL condition. Our analysis also explains why a {\em pre-merging} step is beneficial.  We show  computational experiments in two and three dimensions.
\end{abstract}



\section{Introduction}\label{sec:intro}
State redistribution methods were introduced several years ago \cite{origSRD} as a technique to stabilize finite volume updates on cut cells in an otherwise regular Cartesian mesh. It applies a post-processing step to the small cells that are cut by the solid geometry by   ``merging" only temporarily with one or more neighboring cells in a volume weighted fashion and then reconstructing the solution back on to the Cartesian cut cell mesh. The merging and redistribution steps are done in a way that maintains conservation, and is linearity preserving. (The algorithm is presented for a model problem in \cref{sec:simpleEx} of this paper.) The cut cell  volume fraction threshold of one half was demonstrated in two space dimensions to be stable numerically.

A second paper \cite{GAetal} proposed a generalization of state redistribution (SRD) that removed the abrupt threshold for applying the algorithm. It proposed a variation of SRD that gradually shuts off at one half, so that cells that are only slightly below the threshold receive only a little stabilization. A test case of a supersonic vortex flow with an exact solution showed that this generalization, as one would expect, was less diffusive than the original SRD which did not account for cell volume. Both central merging (a cell combines symmetrically with its left and right neighbors), and normal merging (a cell temporarily combines with the cell in the normal direction) were studied.
The error in central merging was halved with the continuous cut-off. Even    normal merging, which is less diffusive than central merging,  showed a 10-15\% improvement in error. 
That paper also extended SRD to three space dimensions, and to the more complicated equations of low Mach number flow with diffusion and compressible flow with combustion.

Cell merging (see \cite{JiEtAl:2010,saye2017implicit} for some recent references)  is a common alternative  for treating  small cut cells in an embedded boundary mesh. In the most common variation of cell merging, cells are permanently combined into a union of polyhedra with sufficient volume, and the finite volume update is applied to the merged cells. One advantage of cell merging over SRD is that merged cells can be updated in a monotone way, at least for first order linear advection, and as much as possible for second order schemes using appropriate limiters, time steps and time steppers. An advantage of SRD over cell merging is that it avoids the global algorithms that are  applied in the preprocessing step to find good combinations of cells to merge. SRD only temporarily merges cut cells with neighboring cells, in either the direction normal to the wall  (the preferred approach) or in a symmetric fashion called central merging, and it remains stable when two or more cells overlap the same neighboring cell. SRD can be applied in a completely local fashion.

Motivated by the work in \cite{GAetal} we generalize the approach to develop a framework where the weights can be chosen to continuously approach zero when the volume fraction of a cut cell reaches a threshold. The framework also  ensures a monotone update in most situations, like cell merging does.  
We analyze the new weighted algorithm and show that it is both total variation diminishing (TVD) and  monotonicity preserving on a model problem. (Since the scheme is not translation invariant, and has different coefficients on different cells, it is not sufficient to have positive coefficients for monotonicity preservation, nor does TVD imply it in this situation.) The framework for deriving the new weights can also be used to give the original weights, and also those in \cite{GAetal}, giving some insight into when they are not monotone. The analysis also gives some insight into why the original weights work so well. Nevertheless there are reasons one might prefer the monotone scheme, especially since it has negligible additional cost.

The authors are aware only of \textit{h}-box methods \cite{mjb-hel-rjl:hbox}, and Domain of Dependence methods \cite{engwer2020stabilized}, which are conservative, monotone, and TVD.  
There is another cut cell method called flux redistribution \cite{chern:colella} that is also conservative, very simple to implement, but not accurate. It does have the advantage that the amount of post-processing goes to zero when the solution update goes to zero, since it is the cell update that is re-distributed and not the whole state, as in SRD. Flux
stabilization \cite{GOKHALE2018186} is another approach which takes into account the local CFL number
for improved performance at stagnation points, or where the cut cell is large enough
to only require a small amount of stabilization. We do not yet have a way to incorporate these behaviors.

In the next section we illustrate the new weighted algorithm in a simple setting with a model problem. This sets the stage for the general weighted algorithm, where the notation can be a bit unwieldy.  \Cref{sec:cons} introduces the general formulation and proves it is conservative.  \Cref{sec:MLA123} states a number of monotonicity results, with the proofs in the Appendix for the most common cases, and in the Supplementary material for some others, so as not to interrupt the presentation with long stretches of algebra. Computational examples in two and three dimensions are in \cref{sec:compEx}. Conclusions and directions for future research are in \cref{sec:conc}.

\section{A Simple Model Problem with Weights}\label{sec:simpleEx}

\begin{figure}[h!]
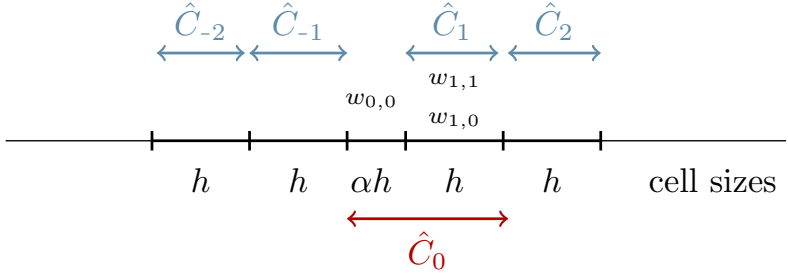

    \centering
    \includestandalone[width=0.8\linewidth]{one_small_right}
    \caption{Model problem with one small cell that merges to the right.}
    \label{fig:1smallCell}
\end{figure}

To illustrate some of the main ideas we use the simple one dimensional model problem of a regular mesh with cells of size $h$ except for one small cell in the middle with size $\alpha  h, \, 0 \le \alpha < 1/2$. This is illustrated in \cref{fig:1smallCell}.  We can take $\alpha < 1/2$ since larger cut cells with $1/2 \le \alpha < 1$ are stable without any post-processing at all. This was proved  in one space dimension in \cite{mjb:stability2} for a model problem with one small cell at an inflow boundary, and has been relied on since then in higher dimensions with no observed stability issues.

For the linear advection equation $u_t + a u_x = 0$, the first order finite volume update  is
\begin{equation}\label{eq:pde}
\widehat{U}_i^n = U^n_i - \frac{ \Delta t} {h_i} \, (f^n_{i+1/2} -f^n_{i-1/2}) =
 U^n_i - \frac{a \Delta t} {h_i} \, (U^n_i -U^n_{i-1}),
\end{equation}
where $f_{i+1/2} = a U_i$.  (We point out that the local Lax Friedrichs (LLF) scheme for linear advection reduces to the upwind scheme for scalar advection. We use LLF for our numerical experiments in \cref{sec:compEx}.)
The numerical update for cell 0 will need to be stabilized.

\vspace*{.15in}
\centerline{\bf Merging Right Example}

We introduce the weighted method in the simplest case by having cell 0 temporarily merge to the right with cell 1.  Other possibilities in one dimension are to merge to the left with  cell -1, or to use central merging which symmetrically merges with both cells -1 and 1.  In \cref{fig:1smallCell} the notation $\widehat C_j$
indicates the merging neighborhoods for each cell $j$. As in original SRD \cite{origSRD} a full cell is its own neighborhood, whereas cell 0's neighborhood, $\widehat C_0 $,  includes both cells 0 and 1. 

The next step is to  define the temporary solution $\widehat Q_j^n $ for each merging neighborhood $\widehat C_j$. Again, all full cells use their own updated solution $\widehat Q_j^n  = \widehat U_j^n $. The solution on $\widehat C_0 $ is a weighted linear combination of $\widehat U_{0}^n $ and $\widehat U_1^n$. The first subscript in the weights refers to the index of the cell contributing its solution, and the second index refers to the cell index of the merging neighborhood being created. 
We define
\begin{equation}\label{eqn:qhat0}
\begin{aligned}
    \widehat Q_0^n  &= \frac
    {w_{0,0} \, h_0 \, }{\widehat V_0 \, h}\widehat U_0^n  + 
    \frac{w_{1,0} \, h_{1} \, } {\widehat V_0  \, h}\widehat U_{1}^n  \\[.08in]
    &= \frac
    {w_{0,0} \, \alpha \,}{\widehat V_0 }  \widehat U_0^n+ \frac{w_{1,0}  \, } {\widehat V_0 }\widehat U_1^n
    \end{aligned}
\end{equation}
    Note that the terms in the numerator are volume weighted, and  $\widehat V_0$ is the correspondingly weighted non-dimensional volume fraction in the denominator,
\begin{equation}\label{eqn:vwidehat0}
    \widehat V_0 \,  h = w_{0,0} \, h_0 + w_{1,0} \,h_1 = ( w_{0,0} \, \alpha + w_{1,0}) \, h
\end{equation}

The final update assigns the merged solutions back to the original grid. Note that we use the same weights $w_{i,j}$ in \eqref{eqn:qhat0} that relate two cells as we do in \eqref{eqn:simple2}, only now the first index in the weights refers to the cell receiving the update, and the second index is the merging neighborhood contributing to it.  This relationship  was also the case for the weights in the original SRD algorithm \cite{origSRD}, and is needed to preserve constants, so we build this in to these new weights too.
The weighted updates are:
\begin{equation}\label{eqn:simple2}
\begin{aligned}
    U^{n+1}_0 &= w_{0,0} \widehat Q_0^n, \\[.05in]
    U^{n+1}_1 &= w_{1,0}  \widehat Q_0^n + w_{1,1} \widehat Q_1^n.
\end{aligned}
\end{equation}
We want the final solution to be a positive  combination of the stabilized neighborhood averages, thus we require $0\leq w_{0,0}, w_{1,1}, w_{1,0} \leq 1$. 

By looking at the conservation requirements the weights can be determined. It is only necessary to look at the two cells treated differently, since all others are conservative due to the base finite volume scheme.  Writing the updates at time $n+1$ in terms of time $n$ gives
\begin{equation}\label{eqn:cons}
\begin{aligned}
\alpha U_0^{n+1} + U_1^{n+1} &=  \alpha \, w_{0,0} \, U_0^n \; + (w_{1,0}+w_{1,1})\,U_1^n \; +  w_{0,0} \, \lambda \, f_{\sm1/2}^n    \\[.08in]
&\,-\,  (w_{1,0}+w_{1,1}) \,\lambda \, f_{3/2}^n  + (w_{1,0}+ w_{1,1}-w_{0,0})\lambda \,f_{1/2}^n.   
\end{aligned}
\end{equation}
For conservation, we must have 
\begin{equation}\label{eqn:wts1}
\begin{aligned}
    w_{0,0} = 1 ,  \\[.05in]
    w_{1,0}+w_{1,1} = 1.
    \end{aligned}
\end{equation}
Then \eqref{eqn:cons} reduces to 
$$
\alpha U_0^{n+1} + U_1^{n+1} = \alpha U_0^{n} + U_1^{n} + f_{3/2}^n-f_{\sm1/2}^n,
$$
showing that this version of  SRD is conservative.

Equation \eqref{eqn:wts1} still has one free parameter, which we can take to be 
$w_{1,1}$. 
By substituting the upwind flux  $f(U_i,U_{i+1}) = U_i$ into the above formulas and defining $\lambda = \frac{a \Delta t}{h}$, we can  write the updates as 
\begin{equation}\label{eqn:simpleMergeRight}
    \begin{aligned}
        U_0^{n+1} &= \frac{
        \lambda \, }{\widehat V_0}U_{\sm1}^n+ \frac{ (\alpha-\lambda w_{1,1})}{\widehat V_0} \, U_0^n + 
        \frac{(1-\lambda)(1-w_{1,1})\,}{\widehat V_0} U_1^n,  \\[.09in]
        U_1^{n+1} &=(1-w_{1,1}) \frac{
         \lambda \,}{\widehat V_0} U_{\sm1}^n + 
        (1-(1-\lambda)w_{1,1} )\frac{\alpha }{\widehat V_0}\, U_0^n   \\[.09in]
        & + \frac{(1-\lambda) (1-w_{1,1}(1-\alpha)) }{\widehat V_0}\, U_1^n.
        \end{aligned}
\end{equation}
We see that all multipliers are positive except potentially $(\alpha-\lambda w_{1,1})/\widehat V_0$. 
By choosing $ 0 \le w_{1,1} \leq \alpha$ this is positive too, since we must have $\lambda \le 1$ for stability. 

More generally, if there are more cut cells and only those with $\alpha<0.5$  are merged, then for monotonicity on the unmerged cut cells we will need to reduce the CFL (or alternatively  merge all cut cells with volume fraction $\le$ CFL).
Remember that the upwind scheme on regular cells of size $h$ has positive coefficients for $\lambda \le 1$, since
$U_i^{n+1} = (1-\lambda) \, U_i^n + \lambda \, U_{i-1}^n$.
If cells with volume fraction greater than one half are not merged, for them to be monotone we need to take $\lambda \leq 0.5$. In this case the coefficient $\alpha-\lambda w_{1,1}$ is positive if 
 $w_{1,1} \leq 2 \alpha$. In either case the update concludes by taking $w_{1,0} = 1 - w_{1,1}$ using \eqref{eqn:wts1}. 
Note that the weights approach zero continuously as $\alpha \rightarrow \alpha_\text{target}$, where the specified target threshold is either 0.5 or 1.0.

If we instead take $w_{1,1}=0$ and $w_{1,0}=1$, the resulting scheme is another way to describe  cell merging. If  $U_0^n = U_1^n$ initially, then the update formulas \eqref{eqn:simpleMergeRight} are identical and give $U_0^{n+1} = U_1^{n+1}$.
In both cases, the update coefficients are all positive and sum to one, so there will be no new extrema during time-stepping. One can also show this scheme is monotonicity preserving and TVD, but this must be shown directly since it does not follow from positive coefficients on an irregular grid.

For comparison, the original SRD algorithm used $w_{1,1} = w_{1,0} = 1/2.$
Equation \eqref{eqn:simpleMergeRight} shows this is not a monotone scheme though we do find it to be well-behaved in practise.
For the generalized $\alpha \beta$ SRD scheme in \cite{GAetal},
the weights depend on the volume fraction threshold $\alpha_\text{target}$ below which SRD is activated.
For $\alpha_\text{target} = 1/2$, the corresponding weights are  $w_{1,1} = 3/4+\alpha/2$ and $w_{1,0} = 1/4-\alpha/2$. This is also not monotone but works well in two and three dimensions.  

\vspace*{.15in}
\centerline{\bf Merging Left Example}
 
The above example of merging right gives some intuition, but does not always reveal when SRD is monotone. For example, merging to the left results in the following one step update analogous to \eqref{eqn:simpleMergeRight},
\begin{equation}\label{eqn:simpleMergeLeft}
    \begin{aligned}
        U_0^{n+1} &= \frac{
         \alpha-\lambda } {\widehat V_0}\, U_0^n + \frac{\lambda + ( 1-w_{\sm1,\sm1})(1-\lambda)} {\widehat V_0} \, U_{\sm1}^n+ 
      (1-w_{\sm1,\sm1}) \frac{  \lambda} {\widehat V_0}\, U_{\sm2}^n  \\[.08in]
        U_{\sm1}^{n+1} &= \frac{       
        (\alpha-\lambda)(1-w_{\sm1,\sm1} )}{\widehat V_0} \, U_0^n +
        \frac{1-w_{\sm1,\sm1}(1- \alpha +\alpha \lambda)}{\widehat V_0} \, U_{\sm1}^n \\
        & \: + (1-w_{\sm1,\sm1}+\alpha w_{\sm1,\sm1})  \frac{\lambda }{\widehat V_0}\, U_{\sm2}^n  .
    \end{aligned}
\end{equation}
The unknown weight now is $w_{\sm1,\sm1}$ and $\widehat V_0 = \alpha+1-w_{\sm1,\sm1}$.  In both equations in \eqref{eqn:simpleMergeLeft} the coefficient in front of the $U_0^n$ terms include $(\alpha-\lambda)$. A time step based on the full CFL step would lead to a negative coefficient for small $\alpha$.

However a small fix to the algorithm does result in a monotone scheme. While looking at $U^{n+1}$ as a function of $U^n$ does not result in a positive scheme, looking at  $\widehat Q^{n}$ as a function of $\widehat Q^{n-1}$ is provably positive. This is stated in more detail in \cref{sec:MLA123} and proved in the Appendix.   
Since $U^{n+1}$  is a convex combination of $\widehat Q^{n}$, this step is positive too.
If the initial data is stabilized first with an application of SRD (also a convex combination), then the resulting algorithm is  monotone for merging left up to the full CFL $\lambda = 1.$ In other words, given  initial conditions $U^{init}$, we apply the SRD algorithm {\em before} a finite volume update is performed  to get the  {\em pre-processed} initial condition $U^0$. This step is called pre-merging, and is defined in \eqref{eqn:premergeDef} in \cref{sec:MLA123}. Pre-merging  was previously observed to be useful in \cite{GAetal}, but we now understand why it is needed.
We will exploit this in the proofs below.

\noindent{\em Remark:}
If the flow were in the negative direction, the behavior of the merging neighborhoods would switch, and the right-merging neighborhood would need the pre-merging step. However pre-merging would  be necessary for systems with characteristics in both directions, so we always use it.

If the original SRD weights of 1/2 are used instead of the new monotone weights for this left-merging neighborhood, and pre-merging is used, it turns out the scheme is monotone for a reduced CFL $\lambda \leq 1/3$.
We remark that the original SRD scheme  is still stable in the GKS sense  for $\lambda \leq 1$ (proved in Appendix \ref{sec:gks}), but it  can have some under- or overshoots.   

\begin{figure}[h!]
  \centering
  \includegraphics[width=0.45\linewidth]{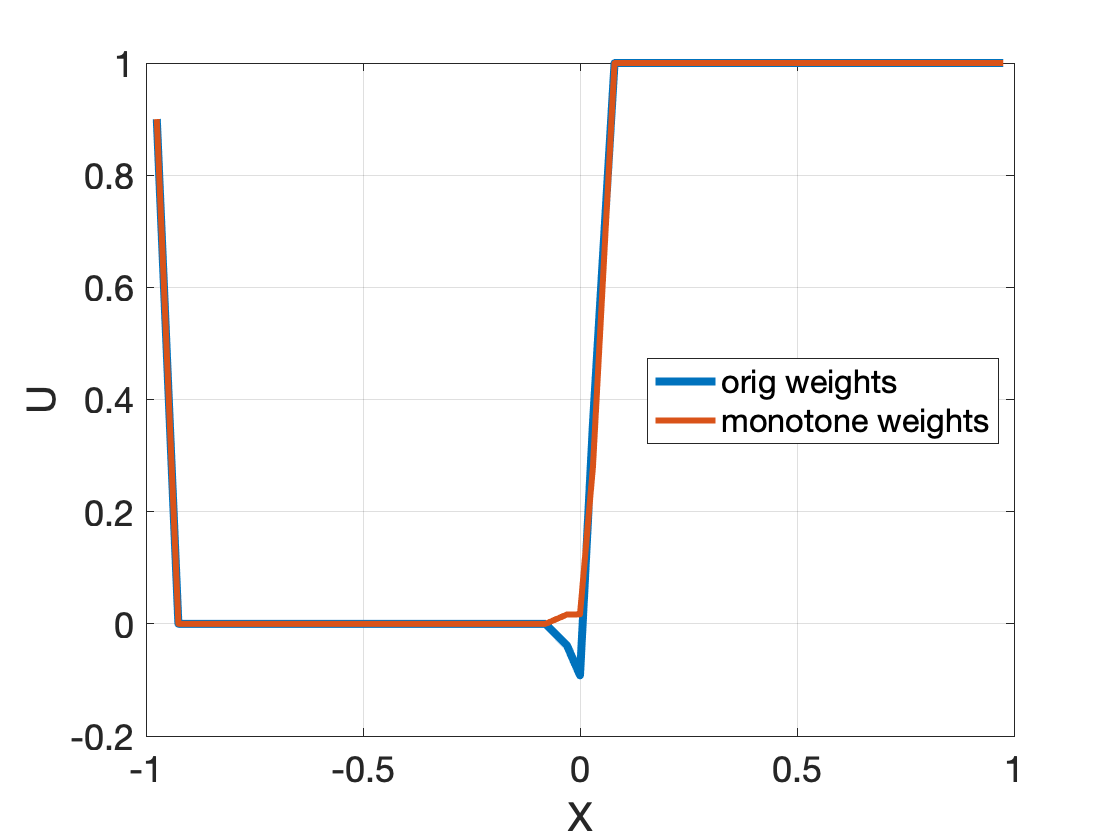}
  \includegraphics[width=0.45\linewidth]{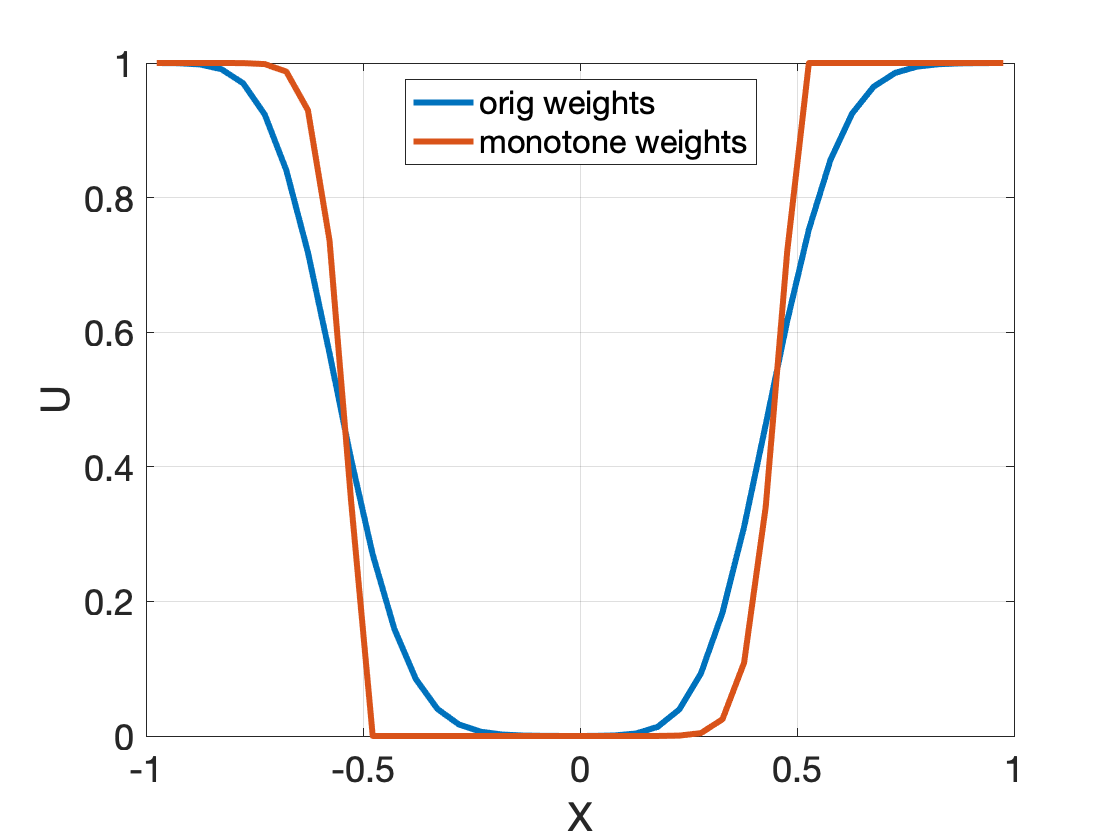}
  \caption{Left figure compares original and new monotone weights after one step using a first-order upwind scheme, with CFL=0.9. The jump is located exactly at the small cell. Right figure compares the monotone weights using the larger time step with many more steps needed for the original weights at the reduced monotone time step, at time $t=0.448$.}
  \label{fig:1dmergeLeft}
\end{figure}

\cref{fig:1dmergeLeft}  shows a square wave advecting to the right on a periodic domain [-1,1]. The initial conditions are
\begin{equation}
    u(x,0) = 
    \begin{cases}
    0 &  x < 0 \\
    1 & x \geq 0
    \end{cases}
\end{equation}
so that the jump occurs exactly at the small cell. The one small cell is located exactly in the middle of the domain. There are   20 full cells on each side, so 40 full cells plus the one small cell discretize the domain.  The small cell has $\alpha = 0.2$. All runs start with a pre-merging step, where the stabilization algorithm is applied to the initial conditions before time stepping begins. The first order upwind scheme \eqref{eq:pde} is used.

The plot on the left shows two curves after one time step, with a CFL of $\lambda=0.9$,   the new monotone weights with $w_{\sm 1, \sm 1}=\alpha$ in red, and the original weights with $w_{\sm 1, \sm 1}=1/2$ in blue.  
The second plot shows the results at time $t=.448$. The 
 original weights now use a CFL of $\lambda = 0.3$ so it is monotone. It takes 30 time steps to reach this time. The new monotone  weights still use a CFL of $\lambda = 0.9$, and only took 10 steps. The original SRD scheme is much more diffusive with this small time step. It  undershoots by -0.09, -0.05,-0.015 for the first three steps. By the $10^\text{th}$ step the diffusive nature of the scheme reduces the undershoot to -0.00015, and would not be visible in these plots.

\section{The Weighted Algorithm and Conservation}\label{sec:cons}

This section defines the new weighted algorithm, gives a condition on the weights,  and proves it is conservative. Conservation holds in 1, 2 or 3 space dimensions by thinking of cell index $i$ as possibly a multi-index.  Since the notation can be confusing, the proof is following by its application to the mesh in  \cref{fig:1smallCell}. 

Recall that
$U_i^{n}$ is a weighted sum of merging neighborhoods $\widehat Q$'s that contribute to cell $i$.  We call this set $W_i$, and use $N_i$ for the cardinality of the set, also referred to as the overlap count for cell i. We write
\begin{equation}\label{eqn:defU}
U_i^{n} = \sum_{j \in W_i} \; w_{i,j} \widehat Q_j^{n-1}.
\end{equation}
For conservation, we will need 
\begin{equation} \label{eqn:sumw}
    \sum_{j \in W_i} \; w_{i,j} = 1.
\end{equation}
There is also a set $M_j$ of cells $i$ such that $\widehat{U}_i$ contributes to $\widehat{Q_j}$,  defined by
 \begin{equation}\label{eqn:defQ}
 \begin{aligned}
 \widehat Q_j^{n-1}  &= \frac{1}{\widehat V_j} \sum_{i\in M_j} \; w_{i,j} \, V_i \, \widehat U_i^{n-1}  , \quad  \text{where} \\
\widehat V_j &= \sum_{i\in M_j} \;  w_{i,j} \, V_i.
\end{aligned}
\end{equation}
The weighted volume is $\widehat{V}_i$. From now on we use $V$ and $\widehat V$ as dimensional variables, with values depending on the number of space dimensions.
 Note that $i \in M_j$ is equivalent to $j \in W_i$, so only one of the sets needs to be formed when implementing the algorithm.  We will prove

{\bf Theorem:} The weighted SRD scheme defined by (3.1)-(3.3) is conservative.

We start with  discrete conservation at time $t^{n}$ and show it reduces to the same quantity at time $n-1$.  
\begin{align}\label{eq:decompose_prev}
\sum_i V_i U^{n}_i &= \sum_i V_i \sum_{j \in W_i} w_{i,j} \widehat Q_j^{n-1},
\end{align}
where $i$ ranges over all cells in the domain.
Rearranging the right-hand-side of the above to group by  $\widehat Q^{n-1}$, and recalling that  $j \in W_i$ is equivalent to $i \in M_j$ we have
\begin{align} \label{eq:decompose}
\sum_i V_i U^{n}_i &= \sum_i \widehat{Q}_i^{n-1} \sum_{j \in M_i} w_{j,i} V_j  .
\end{align}
Recalling the definition of $\widehat V_i$  this simplifies to
\begin{align} \label{eq:intermediate}
\sum_i V_i U^{n}_i  &= \sum_i \widehat V_i \widehat Q_i^{n-1}
\end{align}
Using the definition of the $\widehat Q_i$  in \eqref{eq:intermediate} we obtain
\begin{align} \label{eq:before_decompose}
\sum_i V_i U^{n}_i  &= \sum_i \sum_{j \in M_i} w_{j,i} V_j \widehat U_j^{n-1}.
\end{align}
Regrouping again this time by $\widehat{U}^{n-1}$, we have
\begin{align} \label{eq:intermediate2}
\sum_i V_i U^{n}_i  &= \sum_i V_i\widehat U_i^{n-1} \sum_{j \in W_i} w_{i,j} .
\end{align}
We use the conservation requirement $\sum_{j \in W_i} w_{i,j} = 1$ from \eqref{eqn:sumw}, so \eqref{eq:intermediate2} becomes
\begin{align}
\sum_i V_i U^{n}_i  &= \sum_i V_i  \widehat U_i^{n-1} = \sum_i V_i U_i^{n-1}
\end{align}
since the base finite volume scheme is conservative, thus proving that the overall update with SRD is conservative too. $\blacksquare$

We illustrate the rearranging of terms using the mesh from \cref{fig:1smallCell}.
We have
\begin{equation}
    \begin{aligned}
        M_0 &= \{0,1\}  \hspace{.5in}   W_{0} = \{0\} \\
        M_1 &= \{1\}   \hspace{.65in}     W_{1} = \{0,1\} .
    \end{aligned}
\end{equation}
In words, this is because cell 1 is full, so $\widehat Q_1 = \widehat{U}_1$, and the cell is its own merging neighborhood.  
For cell 1, $U_1^{n+1}$ is a convex combination of $\widehat Q_0^n$ and $\widehat Q_1^n$, the two cell indices that comprise $W_1$.
We obtain  
\begin{equation*}
   \sum_i V_i U^{n}_i = \hdots + V_0 w_{0,0} \widehat Q_0^{n-1}+  V_1 (w_{1,0} \widehat Q_0^{n-1} +   w_{1,1} \widehat Q_1^{n-1}) + \hdots.
\end{equation*}
Grouping by the $\widehat Q^{n-1}$
\begin{equation*}
   \sum_i V_i U^{n}_i = \hdots +( V_0 w_{0,0} +  V_1 w_{1,0} )\widehat Q_0^{n-1} +  V_1  w_{1,1} \widehat Q_1^{n-1} + \hdots.
\end{equation*}
Using the definition of the weighted volume $\widehat V_i$ gives
\begin{equation*}
   \sum_i V_i U^{n}_i = \hdots +\widehat V_0 \widehat Q_0^{n-1} +  \widehat V_1  \widehat Q_1^{n-1} + \hdots,
\end{equation*}
corresponding to \eqref{eq:intermediate}.  Inserting the definition of $\widehat Q$, 
\begin{equation*}
   \sum_i V_i U^{n}_i = \hdots +(w_{0,0}V_0 \widehat U_0^{n-1} + w_{1,0}V_1\widehat U_1^{n-1})  + w_{1,1} V_1 \widehat U_1^{n-1} + \hdots,
\end{equation*}
corresponding to \eqref{eq:before_decompose}.  Grouping by $\widehat U_i^{n-1}$ we obtain
\begin{equation*}
   \sum_i V_i U^{n}_i = \hdots +w_{0,0}V_0 \widehat U_0^{n-1} + (w_{1,0} + w_{1,1})V_1\widehat U_1^{n-1} + \hdots = \sum_i V_i \widehat{U}^{n}_i 
\end{equation*}
corresponding to \eqref{eq:intermediate2}.

\section{Monotonicity for Linear Advection}\label{sec:MLA123}

We start with the general formula for the  weights along with the update formulas from Sec.~\ref{sec:cons} to have them all in one place. First, recall that $\widehat{U}_i$ is the finite volume update on cell $i$. Our SRD results will also use a {\em pre-merging} step, an idea introduced above, where SRD is applied to the initial conditions before time stepping begins.  Referring to \eqref{eqn:defQ}, this step applies SRD to $U^{init}$ instead of $\widehat{U}$  resulting in $U^0$, and we write 
\begin{equation}\label{eqn:premergeDef}
    \begin{aligned}
        U^0 &= SRD(U^{init}) \\
          &=  \frac{1}{\widehat V_j} \sum_{i\in M_j} \; w_{i,j} \, V_i \,  U_i^{init}  . 
    \end{aligned}
\end{equation}

We can now define the weighted SRD algorithm:
\begin{equation}\label{eqn:all}
\boxed{
\begin{aligned}
\widehat Q_j^{n-1}  &= \frac{ 1 }{\widehat V_j}\sum_{i\in M_j} \; w_{i,j} \, V_i \, \widehat U_i^{n-1} \\
\widehat V_j &= \sum_{i\in M_j} \;  w_{i,j} \, V_i,\\
U_i^{n} &= \sum_{j \in W_i} \; w_{i,j} \widehat Q_j^{n-1} \\
w_{i,j} &=   \frac{1}{N_i-1}(1-\alpha_j/\alpha_{\text{target}}) \; \text{  for }  j \in W_{i} \backslash \{ i \} \\[.05in]
w_{i,i} &= 1-\sum_{j \in W_i \backslash \{i\}} w_{i,j}.
\end{aligned}
}
\end{equation}
It is easy to check that the weights sum to 1 as required in \cref{sec:cons}. 
Recall that $\alpha_\text{target}$ is the volume fraction threshold for merging, $N_i $ is the number of neighbors that contribute to the final update of cell $i$, $W_i$ the set of indices of those neighbors with cardinality $|W_i| = N_i$, and $M_j$ is the set of neighbors that contribute to the merging neighborhood of $\widehat{Q}_j$. These are the new weights that give a monotone scheme in the most common  situations.
As with original SRD, the weights can be preprocessed for problems with stationary geometry.

 We will derive monotonicity conditions on the weights 
 for linear advection in two space dimensions, assuming flow parallel to a planar boundary.   This is a useful model for the nonlinear Euler equations with wall boundary conditions, and also for incompressible flow or advection of a tracer.  Of course we won't have this exactly in real cases, but as the mesh is refined the boundary  approaches linearity, and the weighted scheme may still be an improvement over  original SRD. 
Many results apply to both the new weights and the original weights, explaining why they are so well behaved.
There are too many cases to prove monotonicity in 3D even with a planar ramp (except in a very restricted setting), but numerical experiments show that the new weights are well behaved here too.

We have not been able to find a general approach that covers all cases of interest, and our derivations come from looking at particular configurations. In what follows we assume the planar boundary is between 0 and $45^\circ$ to the Cartesian mesh. For larger  ramp angles, the boundary, mesh and cell numberings can all be rotated and reflected into such a configuration. In this section we state the results. 
The derivation of the general formula and its application to the most common case of merging normal to the boundary ($N=2$ overlap count) is relegated to Appendix \ref{sec:mla} so as not to interrupt the presentation of the results with pages of algebra that are not that illuminating. Less common configurations of potential interest (central merging, and an example with an overlap count of 3) are in the Supplementary Material. Our results will show that the full stable time step can be retained for cases with overlap count of 2, the most common case, and is only slightly reduced when the overlap count is 3.

We remind the reader that in 1 space dimension, the upwind scheme has a maximum stability limit $a \, \Delta t /{h} = 1$ and is monotone. If there is one small cell with a volume fraction $\alpha_i$, and no stabilization technique is applied, then for monotonicity the time step must be reduced by $\alpha_i$. The next subsections examine the 2D cases.

\subsection{First order linear advection in 2D}
In this section we state the properties of the first order upwind scheme \eqref{eqn:upwind2d}, where for simplicity we assume the mesh widths are equal in both dimensions.

{\bf Theorem:} When solving $u_t + a u_x + b u_y = 0$ using the first order upwind scheme 
\begin{equation}\label{eqn:upwind2d} 
\begin{aligned}
    u_{i,j}^{n+1} &= u_{i,j}^n - \frac{a \Delta t}{h} \,(u_{i,j} - u_{i-1,j}) -
    \frac{b \Delta t}{h} \,(u_{i,j} - u_{i,j-1}) \\
    &=(1 - \frac{a \Delta t}{h} - \frac{b \Delta t}{h})\, u_{i,j}^n + 
    \frac{a \Delta t}{h} \, u_{i-1,j} + \frac{b \Delta t}{h} \, u_{i,j-1}.
    \end{aligned}
\end{equation}
the following Properties 1-5 hold.

{\bf Property 1:   Base Scheme Monotonicity} 
The first order upwind scheme  is monotone for  
\begin{equation}\label{eqn:cfl1}
\Delta t\, (a/h + b/h) \leq 1.
\end{equation}
This follows if the CFL condition  is met since the update  is a convex combination of solution values with positive coefficients.
See \cite{book:LeVeque02}  for  a discussion of monotonicity.

\vspace*{.15in}

{\bf Property 2:  Unmerged Cut Cell Monotonicity}  
On a Cartesian cut cell grid where cut cells with volume fractions larger than 1/2 are not merged, monotonicity is retained for the full CFL $\Delta t\, (a/h + b/h) \leq 1$.

Property 2 is proved in  \cref{sec:aprop2}. Although the cut cell itself for this case is not merged, its  neighboring cells may be, so the derivation of equation \eqref{eqn:qnp12_appendix} in  \cref{sec:aderiv} still applies. For some intuition, think of the limiting case $b=0$ so the boundary is aligned with the horizontal axis and the flow is parallel to it. The cells are cut only in the $y$ direction, but there is no vertical flow.  The cells are not cut in the horizontal direction, and the full CFL $a \, \Delta t /h \leq 1$ can be taken. Thus the stability limit  is not reduced by the volume fraction as it is in the 1D case. 
If we considered inflow boundary conditions instead of parallel flow there would be restrictions on the CFL. 

\vspace*{.15in}
\centerline{\bf Properties of Weighted SRD  }

The next properties are about SRD. We  examine in detail the results of merging in the direction normal to the embedded boundary.  
We determine positivity by examining the one-step update of the solution averages from one time step to the next, and find conditions that guarantee non-negativity of multipliers in the formula.  In 2D there are only a few  cases to examine.  We also state a few properties for central merging, in the special cases where it is most useful.   The derivation and proofs are in the Appendix and Supplementary Materials, since the algebra is long and not very illuminating.

\begin{figure}[h]
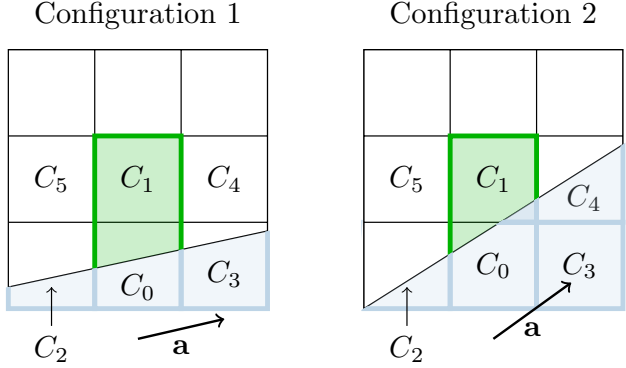

    \centering
    \includestandalone[width=0.3 \textwidth]{2D_cut_square/2D_1_main} \qquad
    \includestandalone[width=0.3 \textwidth]{2D_cut_square/2D_2_main}
    \caption{Merging normal to a planar embedded boundary can result in one of these two possible configurations. In both cases cell 0 has a volume fraction less than 1/2 and an overlap count of 1. Cell 1 has an overlap count of 2 and a volume fraction in larger than 1/2, the target merging volume, so it does not merge with any other cells. }
    \label{fig:2dnormal_A}
\end{figure}

\vspace*{.15in}
{\bf Property 3: Merging Normal to Planar Boundary}
Assume we use merging normal to an embedded boundary with non-negative ramp angle $\le 45^\circ$ (see \cref{fig:2dnormal_A}). Weighted SRD  for configurations 1 and 2, when the overlap count for cell 0 is not more than 2, results in a monotone scheme using the full CFL $\leq$ 1.  The case of overlap count 3 requires a slight reduction of  CFL  $\leq 0.92$, a result obtained with the help of numerics.

The  proof is in \cref{sec:aprop3} for normal merging with overlap count 2, and in the Supplementary Materials for overlap count 3. Substituting  $N=2$,
$\alpha_\text{target} = 1/2$, $W_0 = \{0\}$, $W_1 = \{0,1\}$  into the fully general weights \eqref{eqn:all}, the formulas simplify to
\begin{equation}\label{eqn:weights_normalA}
  \begin{aligned}
w_{0,0} &= 1 \\
w_{1,0} &=  1-2\alpha_0 \\
w_{1,1} &= 1-\sum_{j \in W_1 \backslash \{1\}} w_{1,j} = 1 - w_{1,0} = 2 \alpha_0
\end{aligned}
\end{equation}
\FloatBarrier

\vspace*{.25in}
\centerline{\bf Central Merging, $45^\circ$ Planar Boundary}

The $45^{\circ}$ case is special so we examine it further. It was purposefully used in \cite{GAetal} as a way to test symmetry, which SRD can maintain by using central merging, the subject of the next two results.  The proofs  are in the Supplementary Materials, and use the same approach in the Appendix.

\begin{figure}[h]
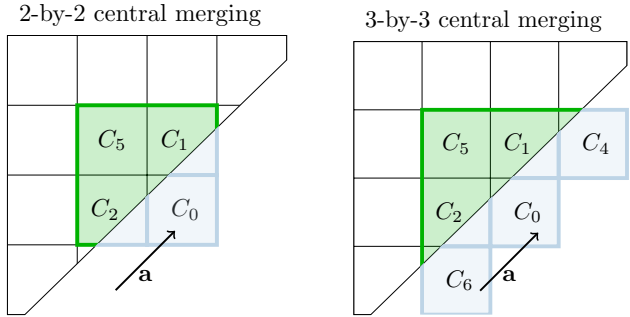

    \centering
    \includestandalone[width=0.3 \textwidth]{2D_cut_square/2D_7_main}\qquad
    \includestandalone[width=0.3 \textwidth]{2D_cut_square/2D_6_main}
    \caption{Two different types of central merging for cell 0 in the $45 ^{\circ}$ ramp case. The first, presented in \cite{GAetal},  uses a smaller merging neighborhood than the 3-by-3 tile on the right.  Left figure has  an overlap count $|W|$ of $N_0=1, N_1 = N_2 = 3, N_5 = 2$, right has counts $N_0=N_1=N_2=3, N_5=2.$ 
    \label{fig:2dcentral}}
\end{figure}

\vspace*{.15in}
{\bf Property 4: 2 by 2 Central Merging} 
For the $45^{\circ}$ ramp angle, where a small cell merges centrally as in  \cref{fig:2dcentral} left (introduced in \cite{GAetal}), the general formula  \eqref{eqn:all} gives a monotone scheme.

Substituting $\alpha_{\text{target}}=1/2$ into the general formula for weights gives
\begin{equation}\label{eqn:weights_central}
  \begin{aligned}
w_{0,0} &= 1 \\
w_{i,0} &=  \frac{1}{N_i-1}(1-2\alpha_0) \;\quad \text{ for } i = 1, 2, 5 \\
w_{i,i} &= 1-\sum_{j \in W_{i} \backslash \{i\}} w_{i,j} \; \text{ for } i = 1, 2, 5. \\
\end{aligned}
\end{equation}
The last equation becomes even simpler for cells 1, 2 and 5:
\begin{equation}
\begin{aligned}
w_{2,2} &= w_{1,1} =  \frac{1}{2} - \alpha_0 \\
w_{5,5} &= 2\alpha_0.
\end{aligned}
\end{equation}
Property 4 is proved in the Supplementary Materials subsection SM2.1.

\vspace*{.15in}
{\bf Property 5:  3 by 3 Central Merging} 
For the $45^{\circ}$ ramp angle, where a small cell merges centrally as in \cref{fig:2dcentral} right,  the weights in \eqref{eqn:all} are not monotone.  Moreover, no choice of weights results in a monotone scheme. 

Property 5 is proved in the Supplementary Materials subsection SM2.2.
For this 3 by 3 central merging case, we do not (yet) have a monotone scheme.  Numerical experiments demonstrate it is  still stable however.

\vspace*{.1in}

To summarize,  we have found weights in the first order case that give positive updates for the solution $U$ at the next time step.  This implies a maximum principle holds, so  global extrema can't grow.  However since the scheme is not translation invariant, and each irregular cell may have a different coefficient, the usual proof of monotonicity preservation and/or a TVD solution does not hold.   
In the Supplementary Material several additional properties are proved. In particular we prove  Property 6, that SRD  in 1D, for the one small cell model problem and merging left with pre-merging, is TVD. This does imply monotonicity preservation or the total variation would increase. 

\vspace*{.25in}

\subsection{Second order linear advection} 

For linear advection  with a second order scheme we use slope reconstruction and a 2 stage TVD Runge-Kutta scheme in time. With appropriate limiting and  a monotone base scheme we find the second order scheme to be well-behaved too, but we have no proofs for this case.

In 1D, on the base grid, von Neumann stability requires a CFL limit $\lambda \le 1$ using central difference slopes, and $\lambda \le 0.5$ for upwind slopes. 
However for both a monotone and TVD scheme the gradient needs to be limited, and the  time step requirement is that the CFL $\leq 0.5$ regardless of how the gradient is initially computed.  
The strong stability preserving (SSP) Runge-Kutta results (see e.g.~\cite{gottlieb2001strong}) state that if a scheme is TVD using forward Euler, then higher order SSP schemes are also TVD since they are a convex combination of forward Euler updates. A second order accurate SSP method is
\begin{equation}\label{eqn:rk2}
    \begin{aligned}
         u^1 &= u^n + \Delta t \, L(u^n) \\
         u^2 &= u^1 + \Delta t \, L(u^1) \\
        u^{n+1} &= \onehalf u^n + \onehalf u^2 
    \end{aligned}
\end{equation}

On the regular cells on the base finite volume scheme we use monotonized central difference gradients. For the irregular cells, in 1D we use minmod and in 2D we use Barth Jespersen as described in \cite{origSRD}. We use the same limiting procedure when reconstructing the gradient on the merging neighborhoods to put the solution back on the Cartesian grid. We always start with the pre-merging step.  
In 1D  we observe numerically that the resulting scheme has no new extrema  with a CFL $\lambda \leq 0.5$ when merging left or right, but have not proved this. The 2 stage scheme has a larger CFL stability limit, but to guarantee that the intermediate stage values are monotone, the forward Euler time step  $\lambda \leq 1/2$ is required and it is  tight. 

As with the original weights in \cite{origSRD} for the second order scheme, the choice of neighborhood averages used to define and limit the slope on a merging neighborhood is crucial. To remind the reader, the second order scheme calculates the centroids of the merging neighborhoods, and a gradient through the cell centroid is evaluated at a neighbor's centroid to calculate its contribution to the final update.   We give an example here using the same model problem as in \cref{sec:simpleEx} with pre-merging and merging left,  but using a second order weighted SRD method and the weights in \eqref{eqn:all}.
We test two different approaches  on neighborhood 0:
\begin{align}
    \widehat{\sigma}_0 &= \text{minmod}\left( \frac{\widehat Q_1-\widehat Q_0}{\hat x_1 - \hat x_0}, \frac{\widehat Q_1 - \widehat Q_{\sm 1}}{\hat x_1 - \hat x_{\sm 1}}, \frac{\widehat Q_0-\widehat Q_{\sm 1}}{\hat x_{0} - \hat x_{\sm 1}} \right),  \label{eq:l1}\\
    \widehat{\sigma}_0 &= \text{minmod}\left( \frac{\widehat Q_1-\widehat Q_0}{\hat x_1 - \hat x_0}, \frac{\widehat Q_1 - \widehat Q_{\sm 1}}{\hat  x_1 - \hat  x_{\sm 1}}, \frac{\widehat Q_0-\widehat Q_{\sm 2}}{\hat x_{0} - \hat x_{\sm 2}} \right), \label{eq:l2}
\end{align}
where the weighted neighborhood centroid is $\hat x_j =  (\sum_{i \in M_j} w_{i,j} V_i x_i)/\hat V_j$ and the physical centroid of the cell is $x_i$.
The second argument in the minmod function is the unlimited reconstructed slopes, and the first and last arguments are computed from forward and backward differences.
Numerical experiments reveal that limiter \eqref{eq:l1} is not monotone, with overshoots of approximately 10\% using a Heaviside function as initial conditions. Limiter \eqref{eq:l2} results in a monotone scheme.  To implement this, in forming the gradient 
we test for weighted centroids that are at least a distance $h/2$ apart or expand the stencil. This is not the case for $\hat x_0-\hat x_{-1} = \mathcal{O}(\alpha h)$.
We also follow this procedure for gradient construction on the base grid when updating the cut cells  (see eq. \eqref{eqn:gradstencil} for an example).
In 2D we apply this distance requirement dimension by dimension.

\section{Numerical experiments}\label{sec:compEx}
We present a convergence study with the new weights in two space dimensions. We then present a shocked flow example in 2D, followed by a smooth  example in 3D with more complex geometry. 
All experiments solve the Euler equations using the local Lax Friedrichs numerical flux. Reflecting boundaries are implemented by evaluating the pressure at the boundary midpoint. For problems with inflow and outflow boundaries, ghost cells are initialized with suitable solution values, either exact inflow or extrapolation outflow.

\subsection{Supersonic Vortex 2D Convergence Test}
First we look at the accuracy of the new weights, with gradients, using an exact solution to the Euler equations in the well-studied  2D supersonic vortex problem (see e.g. \cite{aftosmis:acc}).
We compare the original SRD weights (based only on the number of neighbors) with the new monotone weights.  As with the $\alpha \beta$ weights in \cite{GAetal}, the new weights have less diffusion and reduced error. Here we only study normal merging since that is the preferred approach over central merging. We also use this example to demonstrate the improvement that comes from using a second order accurate gradient at the irregular cells instead of a first order accurate one. Note that the gradient of the full cells in the interior of the domain uses central differences so it is already second order accurate. On the 216 cell mesh the smallest volume fraction was $8.1 \times 10^{\sm6}$. 

\begin{figure}
    \centering
    \includegraphics[width=.9\textwidth]{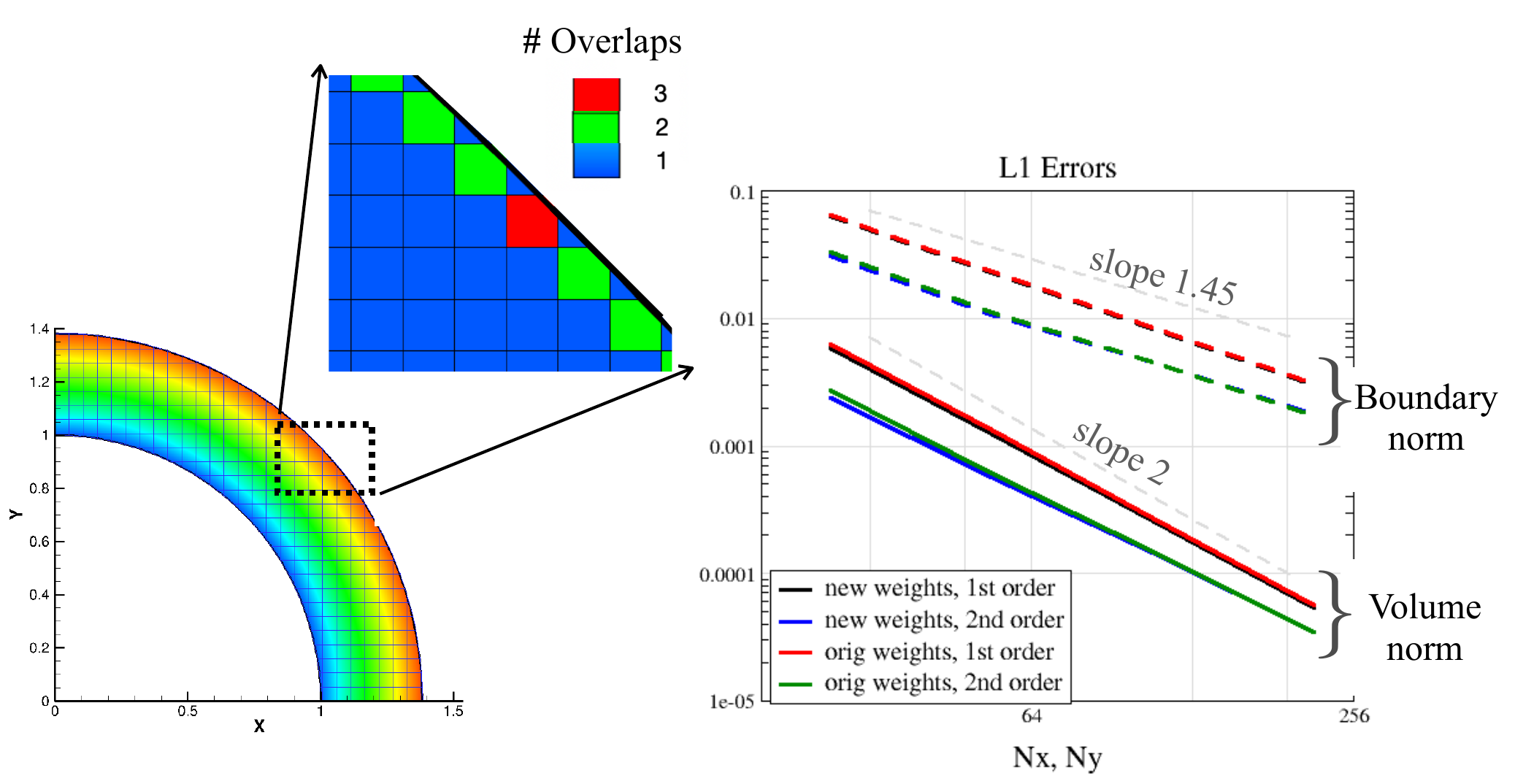}
    \caption{Supersonic vortex test case. Left - density profile of test case. Zoom shows the size of the final updates ($N_i$). On this mesh only one cell has update size 3, with three cells contributing to its final update $U^{n+1}$.  Right - convergence results for normal merging, comparing the new weights and the original weights, with both first order and second order accurate gradient reconstruction. Both the volume norm and boundary norm of the errors are presented. }
    \label{fig:ssvfig}
\end{figure}

The geometry, as in \cite{GAetal}, is depicted in \cref{fig:ssvfig}. With normal merging most cut cells have update size $|W_i|$=2, with the exception of one cell that two cut cells merge with, shown in the zoom. We advance to steady state using  second order Runge-Kutta time-stepping as
in \eqref{eqn:rk2}. The procedure we use to compute first order accurate gradients uses a least squares procedure to fit a linear function for each irregular cell. A cell is irregular if it is cut {\em or} if the cell is full but the central difference stencil has a cut cell.  For the second order gradients, we use a pointwise quadratic reconstruction to fit a quadratic through the irregular cell centroids. The quadratic terms are then dropped and only the gradient is used for reconstruction in time stepping. All gradient reconstructions use primitive variables and no limiting.

For stability, if the gradient stencil doesn't have enough cells it is augmented in one or both directions as follows. 
In any coordinate direction we compute
\begin{equation}\label{eqn:gradstencil}
\max_{i \in \text{stencil}}(i) -\min_{i \in \text{stencil}} (i) \ge 
\begin{cases}
1 \quad  \text{for first order gradients}  \\
2 \quad  \text{for second order gradients}
\end{cases}
\end{equation}
where $i$ is a cell index in the stencil, applied separately in the $x$ and $y$ directions. The stencil is initially of length 3 in each direction, but some of those cells may be in the solid geometry and are not included in the stencil.
This is an easy test to perform on a Cartesian grid. 
In the experiments below, we use the same procedure to compute the limited gradient of the merging neighborhoods as on the base Cartesian grid. This second gradient calculation is used in the final update the of the Cartesian solution from the merging neighborhoods, before taking the next finite volume step.

\cref{fig:ssvfig} shows the L1 norm of the volume error in density  (solid lines) and the L1 boundary norm (dashed lines), defined as the sum of the absolute value of the errors in each cut cell times its boundary segment length. As in \cite{GAetal} the overall convergence rate is second order, and the boundary norms converges at somewhat less than 1.5. The dashed lines show   second order and 1.45 order convergence on the plots. We point out that the convergence of cut cell methods is not smooth, since the grid is highly irregular.  However 
the second order errors are half that of the first order errors, at very little additional cost.  %
In \cite{GAetal} we compared first and second order accurate gradients using central merging. Since the latter is much more diffusive than normal merging, a bigger difference was observed there.

\subsection{Shock diffraction}
\begin{figure}[h!]
    \centering
    \includegraphics[width=.9\textwidth]{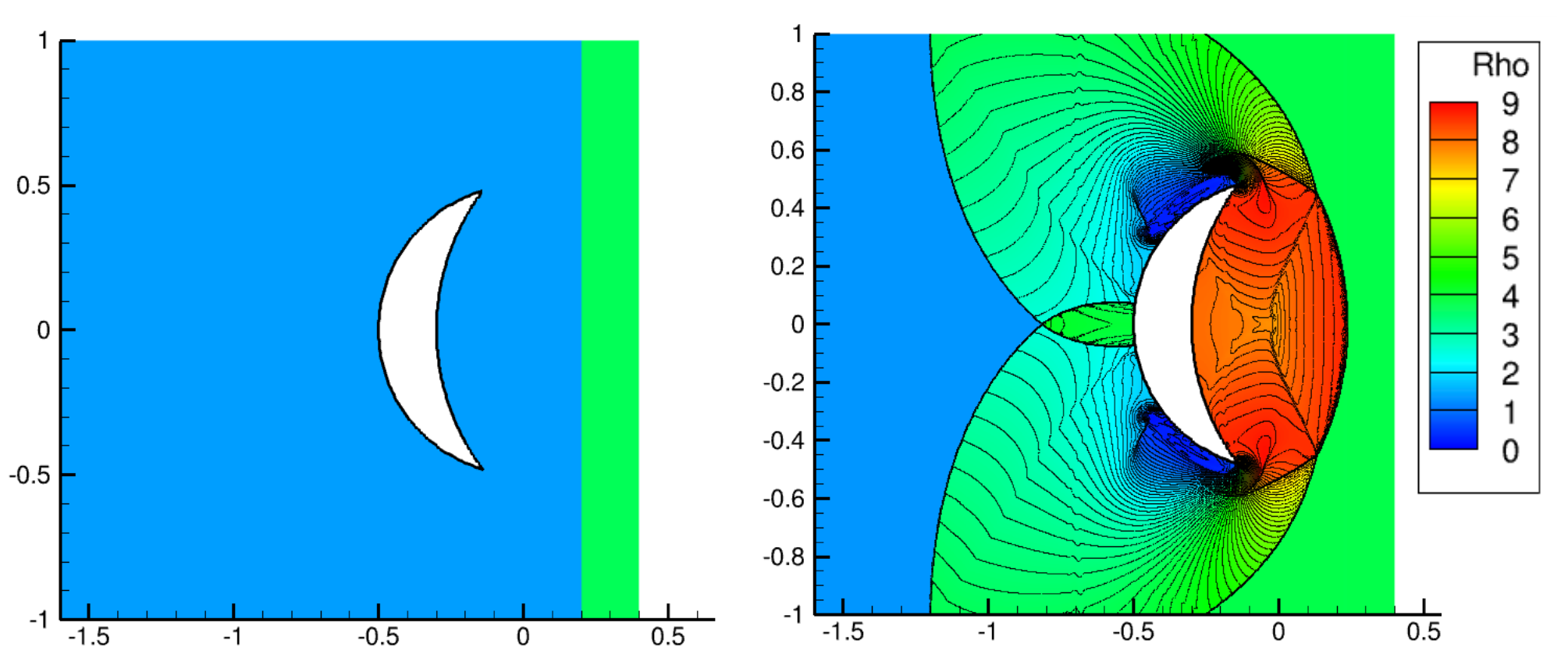}
    \caption{Shock reflection from a crescent. Left figure shows the initial conditions of a Mach 2 shock traveling to the left. Right shows density at time 0.7.}
    \label{fig:crescentVolume}
\end{figure}

Since one of the reasons for monotonicity is to robustly handle shocks, the next example shows shock diffraction around a crescent shape, with both a concave and convex boundary. The domain is $[\sm1.6, 0.4] \times [\sm1,1]$. The mesh spacing is chosen so that the tips of the crescent are located at a gridpoint, so that each refinement still sees the same geometry except for the last cell. Since the geometry is thin at the corners, care must be taken to not include points on the opposite side in the gradient stencils. For this reason we use first order accurate gradients with a 3-by-3  stencil. Since we use normal merging, there is no danger of including cells from the other side in forming neighborhoods. In this example all cut cells only merge with one neighbor, and all final updates include at most one other cell. All the neighborhood  and update sizes are 2 or less in this case. The initial conditions are a Mach 2 shock approaching from the right, as shown in \cref{fig:crescentVolume} left. The solution at time 0.7 is shown on the right. We use the 2 stage Runge-Kutta time stepper with a CFL = 0.5. 

\begin{figure}[h!]
    \centering
\includegraphics[width=.8\textwidth]{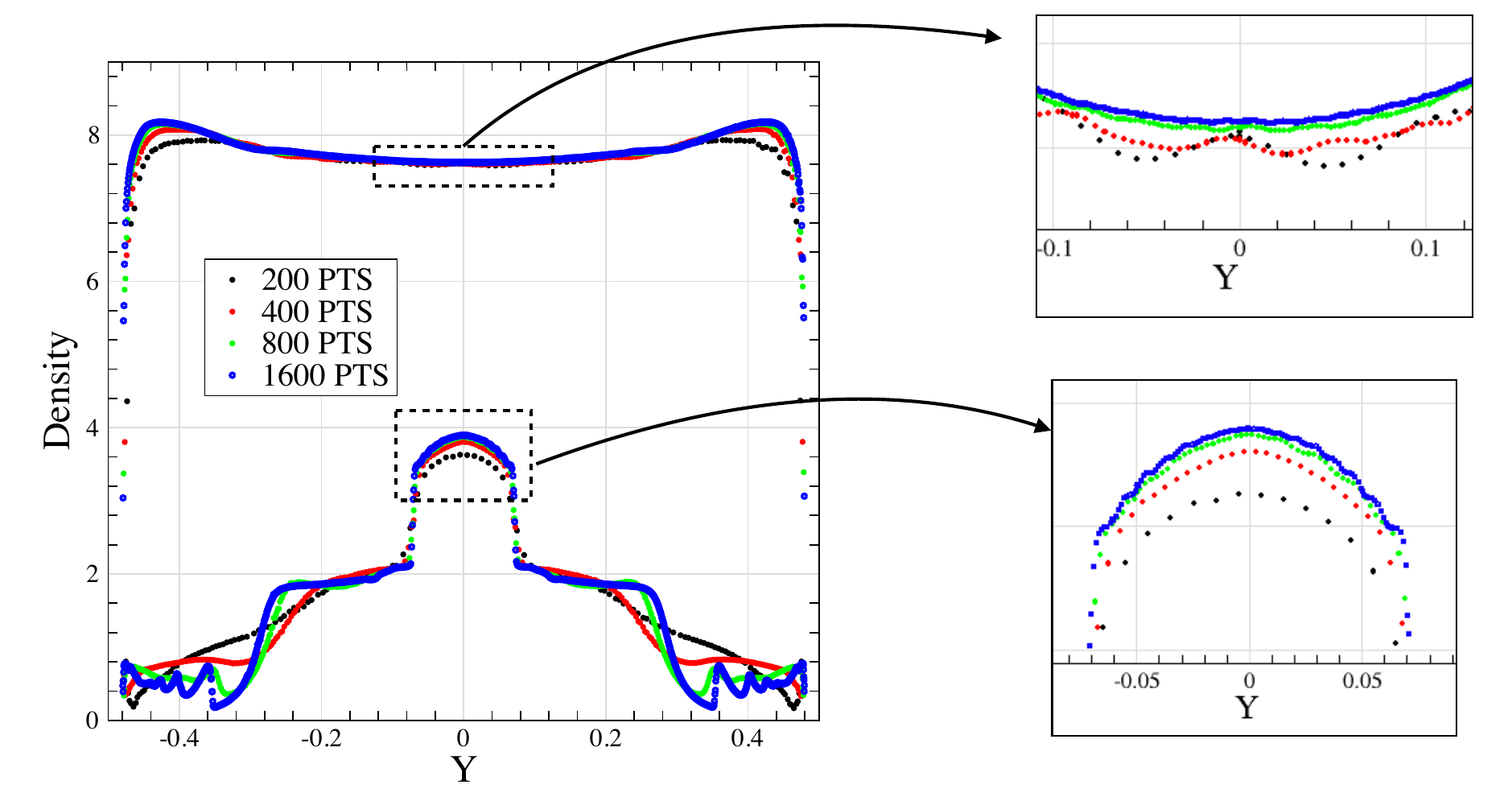}
    \caption{Plot of density along the boundary of the geometry.  The solution at the boundary is reconstructed from the cell centers using first order gradients described in the text. Zooms show the finest two solutions away from the corners are very close.
    }
    \label{fig:crescentBndry}
\end{figure}

A plot showing convergence of density around the surface is in \cref{fig:crescentBndry}, going from 200 to 1600 points in each dimension. The last two data sets, with 800 and 1600 points, are essentially on top of each other away from the singular corner.  The smallest volume fraction in the 1600 cell mesh was    $2.6 \times 10^{\sm17}$, occurring in two cells.

\FloatBarrier

\FloatBarrier

\subsection{Acoustic pulse through trefoil}
In this example, we compute the scattering of an acoustic pulse  by solving the Euler equations inside a trefoil-shaped cavity. The purpose of this example is to show that the weights behave well in 3D curved geometry too, although the analysis in \cref{sec:MLA123} is in two space dimensions and uses planar geometry.  The surface of this complex geometry is meshed with 47,396 flat triangles using Mathematica \cite{Mathematica} from a level set definition \cite{kedia2016weaving}:
\begin{equation}
    \begin{aligned}
        u(x,y,z) &= \frac{2(x + yi)}{1 + x^2+y^2+z^2},\\
        v(x,y,z) &= \frac{2z + i(x^2+y^2+z^2 - 1)}{1 + x^2+y^2+z^2},\\
        \Psi(x,y,z) &= \left |\frac{u(x,y,z)^3}{u(x,y,z)^3 + v(x,y,z)^2} \right|,
    \end{aligned}
\end{equation}
where $i$ is the imaginary number, and $|\cdot|$ is the magnitude of the complex argument.  The trefoil cavity used in this example is given by the implicit equation $\Psi(x,y,z) \geq 1.01$.
\begin{figure}[h]
    \centering
    \includegraphics[height=0.15\textheight]{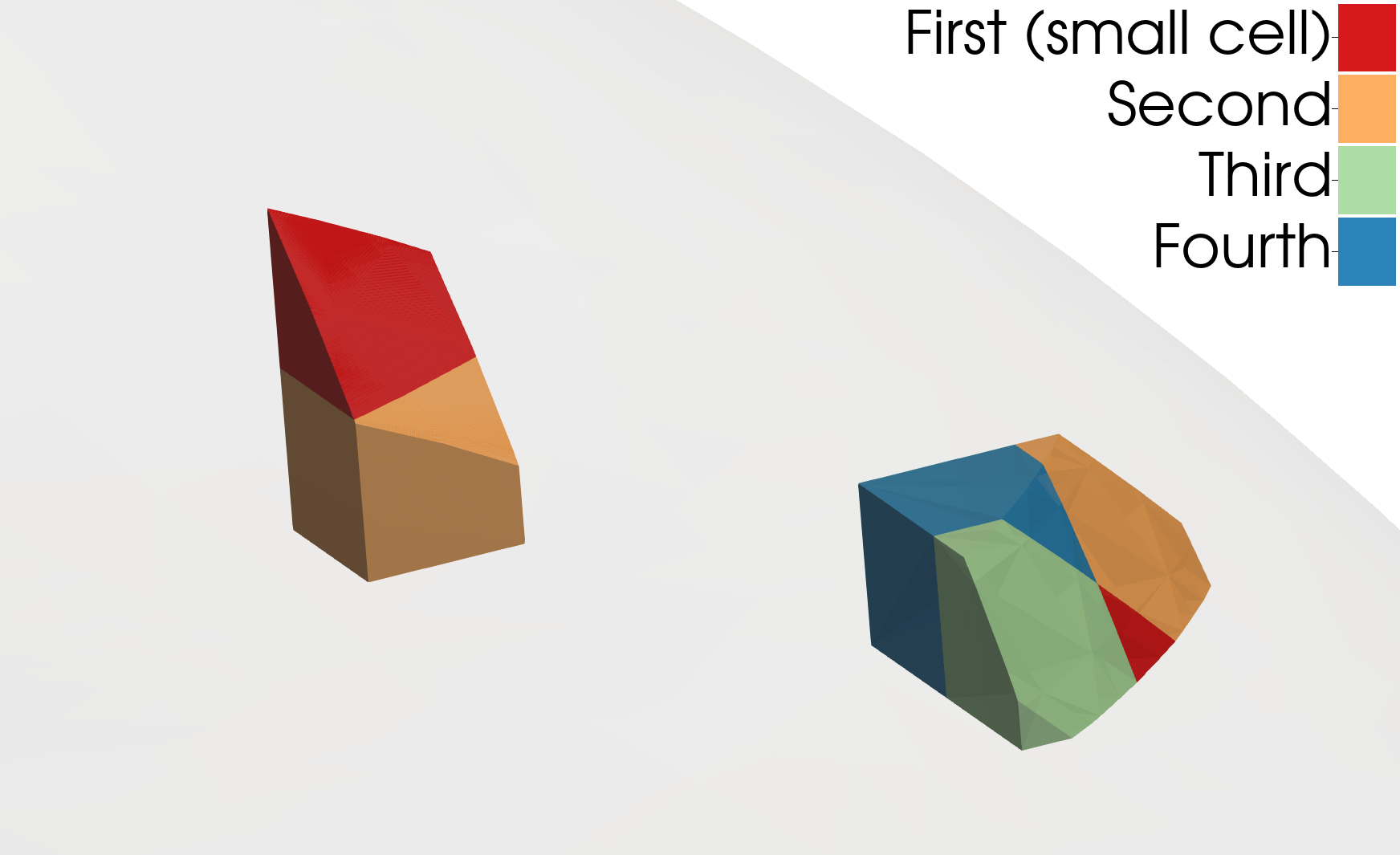} 
    \qquad \qquad \qquad
    \includegraphics[height=0.15\textheight]{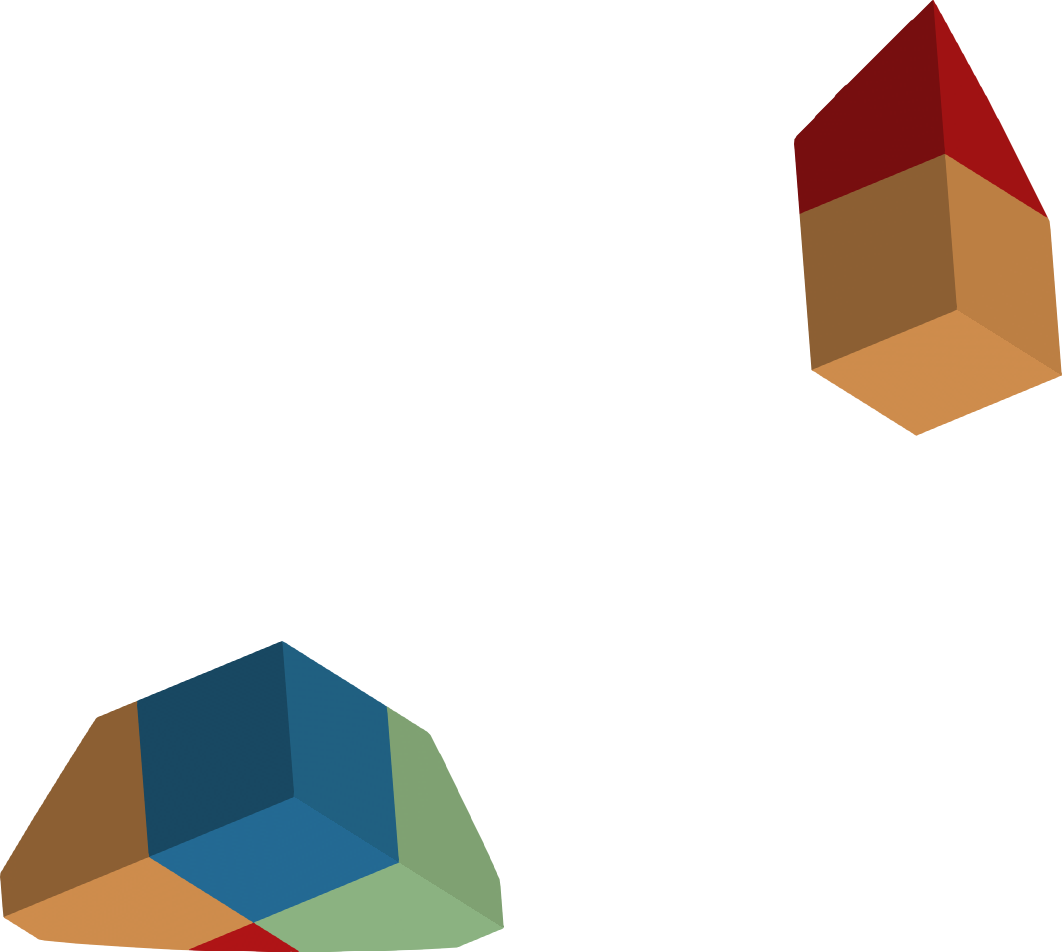}
    \caption{Two possible merging neighborhoods in three dimensions seen from outside the trefoil looking in (left), and vice versa (right).  
    The small cell with volume fraction less than $\alpha_{\text{target}}=1/2$ is colored in red.  
    Subsequent cells added to neighborhood are colored according to the legend.  
    The first neighborhood has only two cells shows standard normal merging.
    The other neighborhood has four cells, since the first and second cells in the neighborhood do not have a large enough volume fraction.
    Thus, the third cell is added, along with the fourth to prevent ``L"-shaped neighborhoods. 
    }
    \label{fig:overlapcounts3D}
\end{figure}
The cut cell mesh is obtained from a $79 \times 79 \times 79$ background Cartesian grid on the domain $[\sm2.25651,3.09301] \times[\sm2.67476, 2.67476] \times[\sm2.67476, 2.67476]$.  
There are 30,874 whole Cartesian cells, 16,983 cut cells in the mesh, and the minimum cut cell volume fraction stabilized by state redistribution is 1.19e-13. 
The complex boundary is first meshed using flat triangles, then this surface mesh is given to the Mandoline mesh generator \cite{Tao:Mandoline:2019} to obtain the three-dimensional embedded boundary grid.

The formula for monotone weights used in the 2D analysis of \cref{sec:MLA123} directly extends to three dimensions and are used here. 
Each merging neighborhood is generated using the normal merging strategy proposed in \cite{GAetal}, with a merging tolerance of $\alpha_{\text{target}}=1/2$.
On this geometry, merging a small cell with a face neighbor that is closest the inward pointing boundary normal will not necessarily result in a neighborhood volume larger than $\alpha_{\text{target}}=1/2$. However, more than 92\% of the time, standard normal was sufficient as shown in \cref{tab:trefoil_overlap_counts} (right).
One solution to this problem is to first merge in the direction of the largest component of the inward pointing normal\footnote{Note, the irregular boundary of a cut cell has many inward pointing normals, so we use a face area weighted average normal.}.
If this neighborhood size is insufficient, then the cell in the direction of the next largest component of the normal is added.  Then, a fourth cell is added to prevent ``L"-shaped neighborhoods, see the supplementary materials or \cref{fig:overlapcounts3D} for an illustration of this kind of neighborhood.  
The overlap counts, neighborhood sizes, and associated frequencies in the cut cell mesh are provided in Table \ref{tab:trefoil_overlap_counts}. The merging strategy used is the reason that we do not observe neighborhoods with 3 cells.
  Most cells have overlap counts of 1 and 2, but the overlap count can be as large as 7.  This is much larger than what we generally observe in two dimensions, where overlap counts usually range between 1 and 3.
Two possible merging neighborhoods are shown in \cref{fig:overlapcounts3D}. 

\begin{table}[]
    \centering
    \begin{tabular}{|c|c|}
    \hline
final update size $|W_i|$ & frequency \\
    \hline
1& 39421 \\
    \hline
2& 6839  \\
    \hline
3& 1263  \\
    \hline
4& 286   \\
    \hline
5& 37    \\
    \hline
6& 6    \\
    \hline
7& 5    \\
    \hline
    \end{tabular}
    \qquad
    \begin{tabular}{|p{35mm}|c|}
    \hline
    \centering merging neighborhood size $|M_i|$ & frequency \\
    \hline
\centering 1& 38738 \\
    \hline
\centering 2& 8463  \\
    \hline
\centering 3& 0  \\
    \hline
\centering 4& 656  \\
    \hline
    \end{tabular}
    \caption{Frequencies of final update size, and neighborhood sizes in the trefoil cut cell mesh, where the background Cartesian grid is $80\times80\times80$.}
    \label{tab:trefoil_overlap_counts}
\end{table}

\begin{figure}
    \centering
    \includegraphics[width=0.3\textwidth]{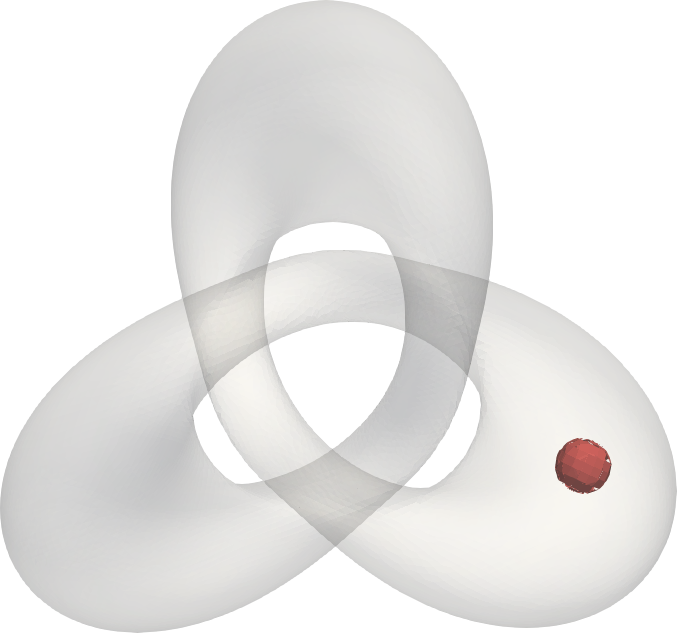}   \quad
    \includegraphics[width=0.3\textwidth]{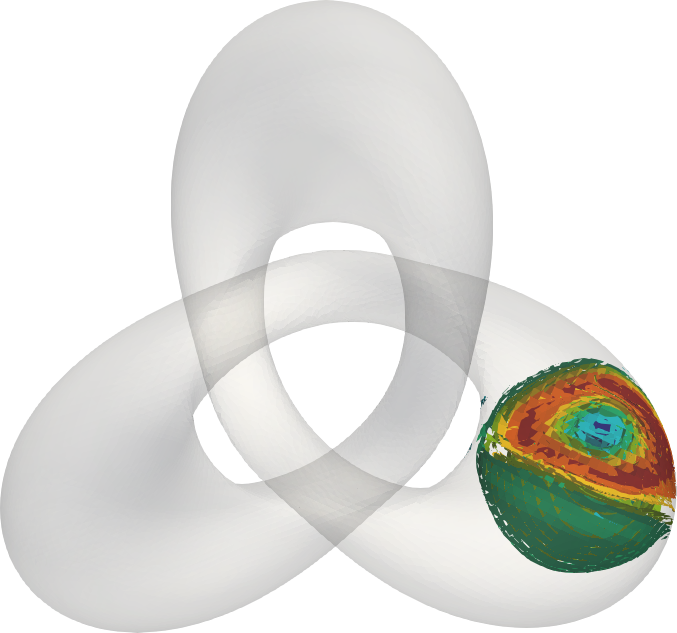}  \quad
    \includegraphics[width=0.3\textwidth]{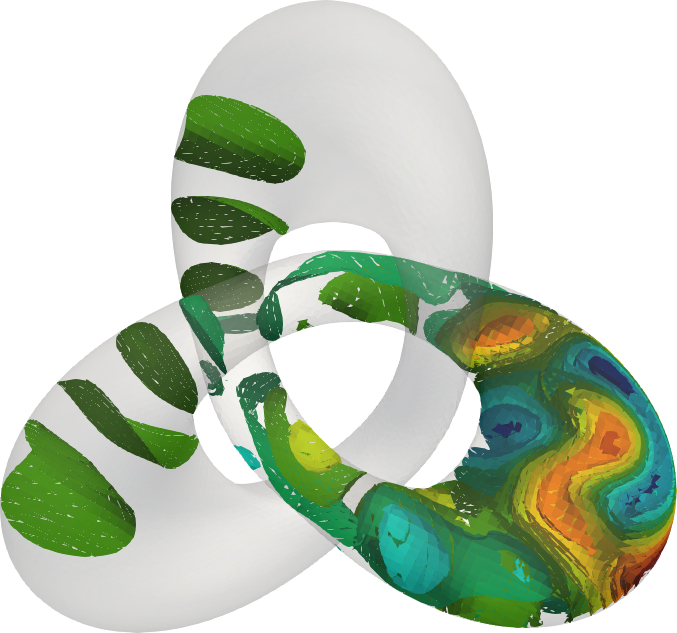}

    \includegraphics[width=0.3\textwidth]{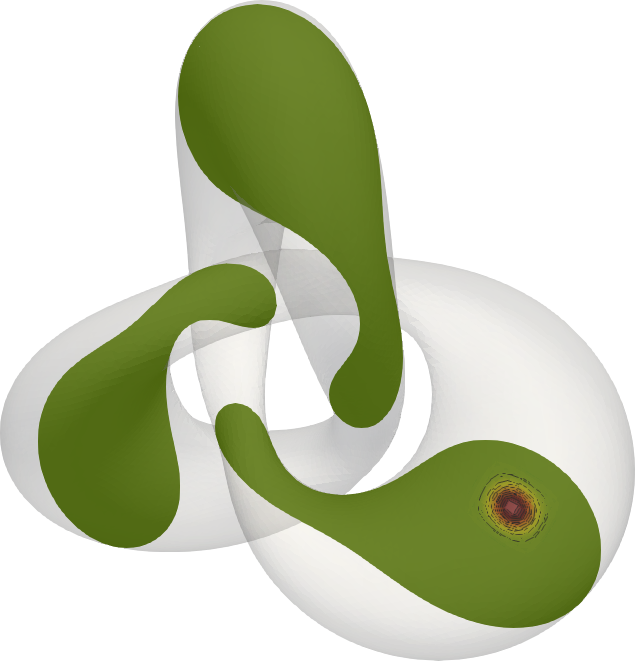}   \quad
    \includegraphics[width=0.3\textwidth]{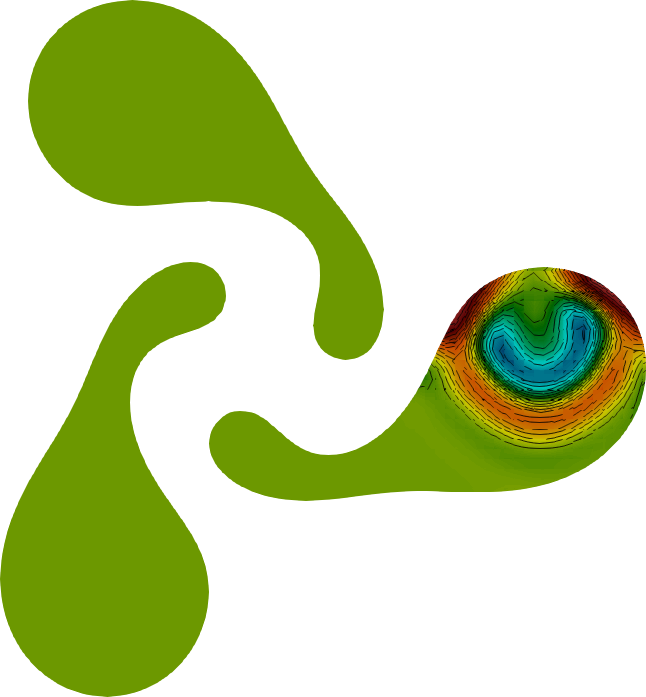}   \quad
    \includegraphics[width=0.3\textwidth]{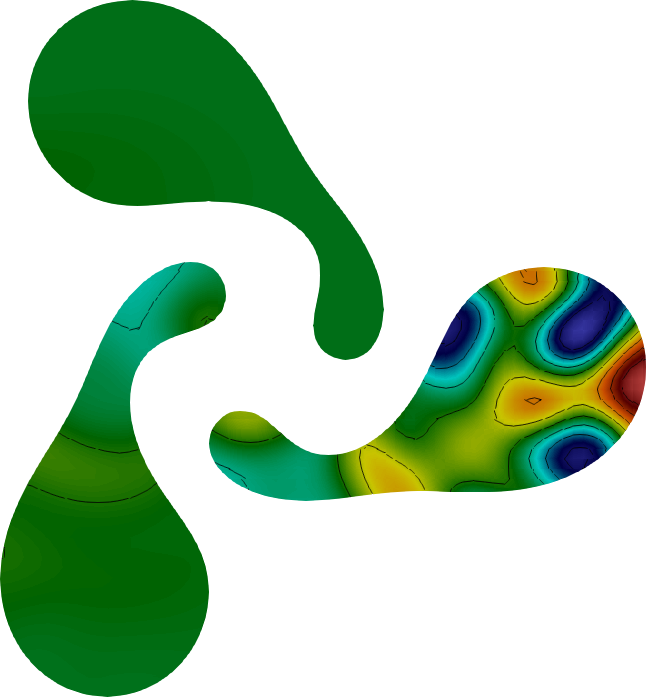}
    \caption{Scattering of acoustic pulse in trefoil cavity.  The trefoil embedded boundary is transparent and encloses the computational domain.  From left to right, each column shows the pressure difference $(P-1/\gamma)/10^{\sm4}$ in the numerical solution at times $t=0, 2\Delta t$ and $22\Delta t$. 
    The top row corresponds to isosurfaces in the pressure difference, while the bottom row shows a slice with isolines on the $z=0.2$ plane.
    The amplitude of the pressure pulse decays quickly, so the color scales in each image differ.
    }
    \label{fig:trefoil}
    \vspace*{-.1in}
\end{figure}

The initial condition is given by
$$
\begin{pmatrix}
\rho\\
u \\
v \\
P
\end{pmatrix} = \begin{pmatrix}
1-1/\gamma + P\\
0\\
0\\
1/\gamma + 10^{\sm4} \exp( - b\left[ (x-x^*)^2 + (y-y^*)^2+ (z-z^*)^2 \right] )
\end{pmatrix},
$$
where $\gamma=1.4$, $b=\text{log}(2)/10^2$, and $(x^*, y^*, z^*)=(1.91867, 0.13616, 0)$ centers the Gaussian perturbation inside the cavity.
No shocks develop since this is a low amplitude perturbation.  Thus, we do not limit and use the full $\text{CFL}=1$.  We reconstruct first order accurate gradients in primitive variables without limiting.
Isobars and slices of the numerical solution are plotted at various times in \cref{fig:trefoil}, where we observe that the solution is smooth up to the embedded boundary.  

\section{Conclusions}\label{sec:conc}
We have proposed a family of weights that generalizes the original State Redistribution algorithm. The original weights were defined using only the number of neighbors of a cut cell. The new weights include the volume fraction,  shut off continuously as the volume approaches the target volume, and are slightly more accurate.  We have analyzed some stability properties of both the new and old SRD weights. The newest formulation is monotone in most situations that occur in practice.  It would be more satisfying to have a general approach that did not rely on specific configurations but we unfortunately did not (yet) find one. 

One SRD property that is still lacking is a way to reduce the amount of post-processing when the solution update goes to zero, as with flux redistribution.
Another useful property might be to  write the weights as functions of the local CFL number $\lambda$ instead of maximizing over it. Finally, we point out that SRD is easily extended to higher order accuracy using a higher-order finite volume scheme. Initial experiments show this is stable.

\appendix
\section{Monotonicity of Linear Advection}\label{sec:mla}
This Appendix first  derives the general formula that covers the configurations in the Theorem of \cref{sec:MLA123}. The second part evaluates the formula in the specific cases to show positivity of coefficients for the stated properties. This includes the cases of cut cells with normal merging (Property 3), as well the large cut cells that are not merged (Property 2).

As with the  1D example in \cref{sec:simpleEx} it is necessary to step from the merging neighborhoods $\widehat Q^{n-1}$ to $\widehat Q^n$ rather than from 
$U^{n-1}$ to $U^n$. The pre-merging step is necessary for monotonicity. It may happen that a particular configuration has positive coefficients without pre-merging but it would not always carry over to parallel flow in the opposite direction without pre-merging.

\subsection{Derivation of $\widehat Q$ Update Formula}\label{sec:aderiv}

First
we show that the neighborhood average $\widehat Q^n_j$ can be written as a convex combination of the merging neighborhoods $\widehat Q^{n-1}_j$ at time $n-1$.
We will look in detail at the contribution of neighborhood average $\widehat Q_j$ on cell $j$ to  cell average $U_i$ for cell $i$, for $j \in W_i$. 
We have  
\begin{equation}\label{eqn:step1}
\begin{aligned}
U_i^n &= \sum_{m \in W_i} \; w_{i,m} \widehat Q_m^{n-1} \\[.08in]
    & =w_{i,j} \, \widehat Q_{j}^{n-1}  + \sum_{m \in W_i,m \neq j} \; w_{i,m} \widehat Q_m^{n-1} \\[.08in]
   &= w_{i,j} \, \widehat Q_{j}^{n-1} \; + \; (1-w_{i,j})\,\widetilde Q_i^{n-1}
\end{aligned}
\end{equation}
where we define
\begin{equation} \label{eqn:qtilde}
        \widetilde Q_i^{n-1} = \frac{1}{1-w_{i,j}} \, \sum_{m \in W_i,m \neq j} \; w_{i,m} \widehat Q_m^{n-1} .
\end{equation}
The $\widetilde Q^{n-1}_i$  are the {\em rest} of the neighborhoods that contribute to cell $i$, excluding cell $j$.
For a concrete example, suppose we have three cells, $C_{\sm 1},C_0$ and $C_1$, with $i=0$, $j=\sm 1$, we can write
\begin{equation}
\begin{aligned}
U^n_0 &= w_{0,\sm 1} \widehat Q^{n-1}_{\sm 1} + w_{0,1} \widehat Q^{n-1}_{1} + w_{0,0} \widehat Q^{n-1}_{0} \\[.08in]
    &= w_{0,\sm 1} \widehat Q^{n-1}_{\sm 1} + (1-w_{0,\sm 1}) 
    \frac{(w_{0,1} \widehat Q^{n-1}_{1} + w_{0,0}\widehat Q^{n-1}_{0})}{(1-w_{0,\sm 1}) } \\[.08in]
    &= w_{0,\sm 1} \widehat Q^{n-1}_{\sm 1} + (1-w_{0,\sm 1}) \widetilde Q^{n-1}_0
\end{aligned}
\end{equation}
Note that $\widetilde Q^{n-1}_{0}$ itself is also a convex combination of neighborhoods and weights that sum to 1, since $w_{0,1}+w_{0,0} = 1-w_{0,\sm 1}$.

Next we go from $U^n$ to $\widehat U^n$, grouping the inflow states for linear advection. 
In 2D (also generalizes to 3D), an update can be written: 
\begin{equation}\label{eqn:inflow}
 \begin{aligned}
     \widehat{U}^{n}_i   &= \left( 1- \frac{\Delta t}{V_i}\sum_{k \in \text{outflow}_i }l_{i,k}h \, \mathbf a \cdot \mathbf n_{i,k} \right) \, U^n_i +\frac{\Delta t}{V_i}\sum_{k \in \text{inflow}_i }l_{i,k}h \, \mathbf a \cdot \mathbf n_{k,i} \, U^n_k \\[.1in]
       &= \left( 1- \frac{\Delta t \, h}{V_i}\sum_{k \in \text{inflow}_i }l_{i,k} |\mathbf a \cdot \mathbf n_{i,k}| \right) \, U^n_i +\frac{\Delta t \,h}{V_i}\sum_{k \in \text{inflow}_i }l_{i,k} |\mathbf a \cdot \mathbf n_{i,k}| \, U^n_k ,
 \end{aligned}
\end{equation}
where $l_{i,k}h$ is the length of the edge between cells $i$ and $k$,  $\text{inflow}_i$ is the set of cell indices that are inflow with respect to cell $i$, and $\mathbf{a}$ is the advection velocity.
Here we have used the divergence theorem $\sum_{k \in \{ cell \, i \, faces\} } \, l_{i,k} \, \mathbf a \cdot \mathbf{n}_{i,k} \, =0$ to write everything in terms of inflow faces and normals and take absolute values. 

Putting $\Delta t$ inside the sums and replacing $V_i = \alpha_i h^2$,  \eqref{eqn:inflow} becomes
\begin{equation}\label{eqn:onestep}
 \begin{aligned}
     \widehat{U}^{n}_i  
       &= \left( 1- \frac{h}{V_i}\sum_{k \in \text{inflow}_i } \Delta t \,  l_{i,k} |\mathbf a \cdot \mathbf n_{i,k}|  \right) \, U^n_i + \frac{ h }{V_i}  \sum_{k \in \text{inflow}_i } \Delta t \,  l_{i,k} |\mathbf a \cdot \mathbf n_{i,k}| \, U^n_k  \\[.05in]
    &= \left( 1- \frac{\lambda_i}{\alpha_i}\right)U^n_i +\sum_{k \in \text{inflow}_i } \frac{\lambda_{i,k}}{\alpha_i} U^n_k .
    \end{aligned}
\end{equation}
where we use the notation   
$\lambda_i = \sum_{k \in \text{inflow}_i}\lambda_{i,k}$ 
where $\lambda_{i,k} = \Delta t ~ l_{i,k} |\mathbf{a} \cdot \mathbf n_{i,k}|/h$. 

The final step to go from the finite volume update $\widehat U^n$ to the merging neighborhoods $\widehat Q^n$. After multiplying both sides of \eqref{eqn:all} by $\widehat V_j/h^2$, and substituting $V_i = \alpha_i h^2$, (in 2 dimensions),  we get 
\begin{equation} \label{eqn:qnp1}
\frac{\widehat V_{j}}{h^2} \widehat Q^{n}_j =  \sum_{i \in M_j} w_{i,j} \alpha_i \widehat{U}^{n}_i
\end{equation}
After some routine algebra pre-multiplying $\widehat{U}^{n}_i$ in \eqref{eqn:onestep} by $w_{i,j}\alpha_i$,  substituting into \eqref{eqn:qnp1}, and using \eqref{eqn:step1}  to replace the ${U^n}$ terms, we finally get the expression for the $\widehat{Q}$ updates in terms of those from the previous step,

\begin{equation}\label{eqn:qnp12_appendix}
\boxed{
  \begin{aligned}
\frac{\widehat V_{j}}{h^2} \widehat Q^{n}_j &=  \sum_{i \in M_j} \biggl[ \left( w_{i,j}^2 (\alpha_i -   \lambda_i) + w_{i,j}\sum_{k \in \text{inflow}_i}w_{k,j} \lambda_{i,k} \right) \widehat Q^{n-1}_j\\
&+\, (w_{i,j}-w_{i,j}^2)  (\alpha_i- \lambda_i) \; \widetilde Q^{n-1}_{i}\\
&+ \, w_{i,j} \sum_{k \in \text{inflow}_i} \lambda_{i,k}\left(1 -w_{k,j}   \right) \,\widetilde Q^{n-1}_{k} \biggr]  
  \end{aligned}
  }
\end{equation}

This is the fully general formula that holds in 1, 2 and 3 dimensions with appropriate definitions of areas and volumes.

\subsection{Monotonicity Proofs}
Next, \eqref{eqn:qnp12_appendix} is adapted to  first prove Property 3 and then Property 2. We examine in detail the possible configurations to show the updates have positive coefficients.

\vspace*{.1in}

\subsubsection{Proof of Property 3}\label{sec:aprop3}

\begin{figure}[h!]
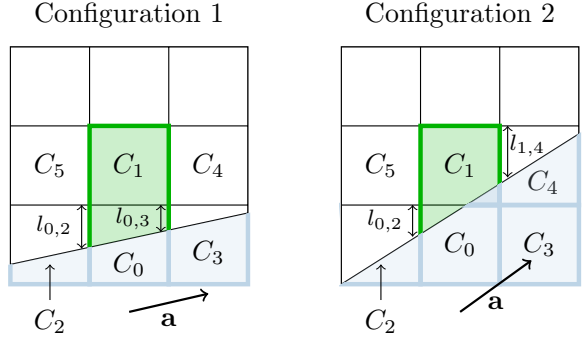

    \centering
    \includestandalone[width=0.275 \textwidth]{2D_cut_square/2D_1} \qquad
    \includestandalone[width=0.275 \textwidth]{2D_cut_square/2D_2}
    \caption{Notation for the two possible configurations in 2D with SRD using normal merging for cell $C_0$ and flow parallel to a planar boundary.}
    \label{fig:2dnormal_appendix}
\end{figure}

In 2D there are two cases for normal merging, illustrated in \cref{fig:2dnormal_appendix} (same as \cref{fig:2dnormal_A} but with more notation).  In both cases the merging neighborhoods contain only two cells from the base grid, and the formulas greatly simplify.

Let  cell $0$ be the index of  the small cell and cell $1$ is the large cell that it merges with.  
In both cases we must have $w_{0,0} = 1$, cell 0 is only in one neighborhood - its own.  The inflow cell for cell 0 is only cell 2. The inflow cells for cell 1 are cell 0 and 5. 
Looking at cell $C_1$ in configuration 2,  $M_1 = \{1\}$ since it is a large cell that doesn't need to merge with other cells. 
The coefficient $w_{5,1} = 0$ since those two cells do not share a neighborhood. $w_{0,1}=0$ because cell 1 does not contribute to cell 0's final update.

For these two cells, \eqref{eqn:qnp12_appendix} becomes 
\begin{equation}\label{eqn:A8}
\begin{aligned}
\frac{\widehat V_{0}}{h^2} \widehat Q^{n}_0  =&   \left(\alpha_0 -    \lambda_0  + w_{1,0}^2 (\alpha_1 -  \lambda_1) + w_{1,0} \lambda_{1,0} \right) \widehat Q^{n-1}_0 \,+\\
& \, (w_{1,0} -w_{1,0}^2)(\alpha_1 - \lambda_1) \widehat Q^{n-1}_{1} 
+  \lambda_{0,2} \widetilde  Q^{n-1}_{2} +w_{1,0} \lambda_{1,5} \,\widetilde Q^{n-1}_{5}  \\[.08in]
\frac{\widehat V_{1}}{h^2} \widehat Q^{n}_1  &=   w_{1,1}^2 (\alpha_1 -   \lambda_1)   \widehat Q^{n-1}_1 +\, \left[(w_{1,1}-w_{1,1}^2)  (\alpha_1- \lambda_1) +  w_{1,1} \lambda_{1,0}\right] \; \widehat Q^{n-1}_{0} + \\
&   +w_{1,1} \lambda_{1,5} \,\widetilde Q^{n-1}_{5}
\end{aligned}
\end{equation}
In both configurations cell 1 has not been merged since it has volume fraction $\alpha_1 \ge 0.5$.
Since there is only one inflow cell face for cell 1 that is different from cell 0,  the definition of $\lambda_{1}$ simplifies to $\lambda_{1,5} =\Delta t \, l_{1,5} a /h$, where the length $l_{1,5} \le 1 $ and $a$ is the $x$-component of the velocity $\mathbf a$. The term $\widetilde Q_5^{n-1}$ represents the neighborhood contributing to cell 5, which might include cell 2 along with cell 5 itself. However since the coefficient in front of this term is positive it does not play a role in determining the monotonicity of the update.
It will be shown below that the term  $\alpha_1-\lambda_1$ in both equations is positive for both configurations if you take a stable time step on the full Cartesian cells in the base grid satisfying \eqref{eqn:cfl1}. 

Thus in \eqref{eqn:A8} the only possibly negative coefficient is  that of $\widehat{Q}_0^{n-1}$ in the first equation. 
 We must find a weight $w_{1,0}$ such that
\begin{equation}\label{eq:monotone_req}
 \begin{aligned}
 \alpha_0 -    \lambda_0+w_{1,0}^2 (\alpha_1 -  \lambda_1) +  w_{1,0}\ \lambda_{1,0}  \geq 0. 
 \end{aligned}
\end{equation}
For flow down the ramp with velocity $-\mathbf a$,  the same inequality results.

\vspace*{.1in}
{\bf Cut Cell Configuration 1:}

We can take the velocity to be 
$\mathbf a = [1, l_{0,2} - l_{0,3}]$, which is parallel to the embedded boundary.  Since $\Delta t$ satisfies 
$\Delta t \cdot (1/h + (l_{0,2}-l_{0,3})/h) = \text{CFL} \leq 1$, the velocity magnitude does not affect this analysis.
For this case we have the volume fraction $\alpha_1 = 1$ (the cell is not cut), and  $\lambda_1 = \frac{\Delta t}{h} \cdot  (1 + (l_{0,2}-l_{0,3}) )$. This gives 
 $\alpha_1-\lambda_1 \geq 0$, and is zero if the maximum time step is taken on the base grid.  Thus it is sufficient that $w_{1,0}$ satisfy
    $$
    \begin{aligned}
 \alpha_0 -    \lambda_0 + w_{1,0} \lambda_{1,0}  \geq 0.
    \end{aligned}
    $$
Solving for the weight, we obtain
\begin{equation}\label{eqn:ineq}
\frac{\lambda_0-\alpha_0}{\lambda_{1,0}}\leq w_{1,0}.
\end{equation}
The left hand side of this inequality  depends on the geometry and flow angle of the channel. We want this estimate to hold for all flow angles less than $45^{\circ}$,  which we can exactly calculate.  Referencing \cref{fig:2dnormal_appendix}, we define the  left and right vertical edges of cell 0 to have edge fraction $0 \leq l_{0,2}, l_{0,3} \leq 1$, $l_{0,2} \geq l_{0,3}$.     The linear advection problem is now completely defined.
We can then evaluate
\begin{equation}
\begin{aligned}
\alpha_0 &= \frac{1}{2}(l_{0,2} + l_{0,3})  \\
\lambda_0 &= \Delta t\frac{l_{0,2}}{h}, \; \; \lambda_{1, 0} = \Delta t \frac{l_{0,2}-l_{0,3}}{h} 
\end{aligned}
\end{equation}
It is easier to analyze using the maximum stable time step, which we do next.
Substituting everything into the left hand side of the inequality \eqref{eqn:ineq},  we obtain
\begin{equation}\label{eqn:ineq2}
\begin{aligned}
\frac{\lambda_0-\alpha_0}{\lambda_{1,0}} &= \frac{ l_{0,2} - l_{0,3} - (l_{0,2}^2 - l_{0,3}^2)}{  2 \, (  l_{0,2} - l_{0,3} )} 
= \frac{1}{2} - \frac{1}{2}(l_{0,2} + l_{0,3}) \\
&= \frac{1}{2}-\alpha_0 \; \leq w_{1,0} \; \leq \, 1.
\end{aligned}
\end{equation}
We have also verified by differentiating with respect to the CFL number that this maximizes the left hand side.

Summarizing, this shows that for this case where $\alpha_0 < 0.5$ in the overlap 2 case,  both  the original SRD weights, where $w_{1,0} = 1/2$, and the new monotone weights from \eqref{eqn:all}, with $w_{1,1} = 2\alpha_0, w_{1,0} = 1 - 2\alpha_0$, satisfy inequality \eqref{eqn:ineq2}.  This helps explain the good performance of the original SRD algorithm in this common configuration. For the previously mentioned $\alpha\beta$ scheme, configuration 1 results in weight $w_{1,0} = \frac{1}{2} (\frac{1}{2} - \alpha_0)$, which does not satisfy \eqref{eqn:ineq2}, but may satisfy \eqref{eq:monotone_req} depending on the time step.

\vspace*{.1in}
{\bf Cut Cell Configuration 2:}

The advection velocity using notation for this configuration is $\mathbf a = [1, l_{0,2}+1-l_{1,4}]$ and $0 \leq  l_{0,2}, l_{1,4} \leq 1$.  Requiring ramp angles less than $45^{\circ}$ implies $l_{0,2}+1-l_{1,4} \leq 1$.  The parameters are:
\begin{equation}\label{eqn:case2}
\begin{aligned}
\alpha_0 &= \frac{1}{2} \, \frac{l_{0,2}^2}{l_{0,2} +1-l_{1,4}}  \; \text{ and } \; \alpha_1 =  1-\frac{1}{2} \, \frac{(1-l_{1,4})^2}{l_{0,2}+1-l_{1,4}} \\
\lambda_0 &= \Delta t \frac{l_{0,2}}{h}, \; \: \lambda_{1, 0} = \Delta t \frac{l_{0,2}}{h} \; \text{  and  } \: \lambda_1 = \Delta t \left(\frac{1}{h} + \frac{l_{0,2}}{h} \right) \\
\Delta t &\left( \frac{1}{h}+\frac{l_{0,2} + 1-l_{1,4} }{h} \right) = \text{CFL} \in (0, 1]
\end{aligned}
\end{equation}

First we show that the multiplier $(\alpha_1 -  \lambda_1)$ of the quadratic term in \eqref{eq:monotone_req} is a positive number.
Unlike configuration 1, it is zero only in special cases such as a 45$^{\circ}$ channel where both $l_{0,2}$ and $l_{1,4} = 0$,  or a horizontal wall.
To simplify notation we define $s = (1-l_{1,4})$ in the algebra below. We want to show
\begin{equation}\label{eqn:A14}
    \alpha_1-\lambda_1 = 1 - \frac{1}{2}\frac{s^2}{l_{02}+s} - \frac{1+l_{02}}{1+l_{02}+s} {\geq} 0 .
\end{equation}
After clearing denominators and collecting terms \eqref{eqn:A14}  becomes equivalent to showing 
\begin{equation}
\begin{aligned}
    2 \, s \, l_{02} + s^2 -s^2 \, l_{02}-s^3  & \geq 0   \qquad \text{or} \\
    s \, (s+2\,l_{02}) - s^2 (s+l_{02}) & \geq 0 
\end{aligned}
\end{equation}
and finally
\begin{equation}\label{eqn:A16}
    \frac{s}{s^2}  \geq  \frac{s+l_{02}}{s+2 \, l_{02}}
\end{equation}
This is clearly true since $s \le 1$ implies the left hand side is $ \ge 1$, and $l_{02}$ and $s$  positive means the right hand side is $ \le  1$. Thus the original expression \eqref{eqn:A14} is non-negative.

It remains to show that the quadratic  \eqref{eq:monotone_req}, repeated here for convenience, is positive,
\begin{equation*}%
 \begin{aligned}
 \alpha_0 -    \lambda_0+w_{1,0}^2 (\alpha_1 -  \lambda_1) +  w_{1,0}\ \lambda_{1,0}  \geq 0. 
 \end{aligned}
\end{equation*}
Here we look at the simpler case where  the overlap count for cell 1 is 2. We can drop the quadratic term with the positive coefficient $\alpha_1 - \lambda_1$  to get the expression   
$$
\frac{\lambda_0-\alpha_0}{\lambda_{1,0}} \leq w_{1,0} 
$$
Using the parameters in \eqref{eqn:case2} the lower bound on $w_{1,0}$ simplifies,
\[
\frac{\lambda_0-\alpha_0}{\lambda_{1,0}} = 1-\alpha_0 - \frac{\alpha_0(2-l_{1,4})}{l_{0,2}},
\]

Remembering that $w_{1,0} = 1-2 \alpha_0$ in the $N_1=2$ case and simplifying further, we need to check that 
 \begin{equation}
 1-\alpha_0 - \frac{\alpha_0(2-l_{1,4})}{l_{0,2}} \leq 1-2\alpha_0.
 \end{equation}
 and it  reduces to requiring $l_{0,2} + l_{1,4} \leq 2$, which is true for any ramp with an angle $\le 45^{\circ}$. $\blacksquare$
 
With the original weights, the same quadratic leads to a slight reduction in the CFL limit to 0.94 for monotonicity.
The Supplementary Material looks at an  artificially constructed case for a planar ramp where we force $N_1=3$ so $w_{1,0} = 1/2-\alpha_0$. In this case numerics show that the CFL limit should be reduced to 0.92 to preserve monotonicity. If the original SRD weights are used, an even smaller CFL of 0.71 is needed.  This helps explain the good performance of the original SRD algorithm in these common configurations.
If a full $\text{CFL} = 1$ is taken, original SRD has negative coefficients, and monotonicity cannot be guaranteed, although it takes a well constructed example to observe it.\footnote{One such example is a 45 degree channel with edge fraction $l_{0,2} = 1/6$, with $U=1$ on one small cut cell and zero elsewhere. There will be a tiny undershoot for a few steps.}

\subsubsection{Proof of Property 2}\label{sec:aprop2}
There are also two possible configurations for large, unmerged cut cells to examine for positivity.

\vspace*{.1in}

{\bf Unmerged Cut Cell $C_1$}

First we look at cell $C_1$ in Configuration 2, already illustrated in \cref{fig:2dnormal_appendix}.
Looking at the cell 1 update in \eqref{eqn:A8} and substituting the parameters in \eqref{eqn:case2} shows the only term that needs to be considered has a coefficient $(\alpha_1 - \lambda_1)$. This was  proved positive in \eqref{eqn:A14}--\eqref{eqn:A16}  above.

\vspace*{.1in}
{\bf Unmerged (Larger) Cut Cell $C_0$}

The other possibility for  a large cut but unmerged cell has the embedded boundary enter from the left face and exit to the right. This would be the case for example for cell $C_0$ in Configuration 1 if the embedded boundary were slightly lower in the $y$ direction so the volume fraction  $\alpha_0 \ge 0.5$. We will use the notation shown in \cref{fig:2dnormal_appendix}.

Simplifying the general derivation \eqref{eqn:qnp12_appendix} for cell 0,  with $M_0 = \{0\}$, inflow$_0 = \{2\}$, $w_{0,0} = 1$, $w_{2,0} = 0$, we get
\begin{equation}\label{eqn:A19}
\frac{\widehat{V}_0}{h^2} \, \widehat{Q}_0^n = (\alpha_0 - \lambda_0) \,\widehat{Q}_0^{n-1} + \lambda_{02} \, \tilde{Q}_2^{n-1}.
\end{equation}
Since $\lambda_{02} = {\Delta t \, l_{02}}/{h}$
is positive, it remains to show that $(\alpha_0 - \lambda_0)$ is.

Substitute the following into  \eqref{eqn:A16} 
\begin{equation}
\begin{aligned}
\alpha_0 &= \frac{1}{2}(l_{0,2} + l_{0,3})\\
   \lambda_0 &= l_{0,2}\frac{\Delta t}{h} \\
 \Delta t & \left(\frac{1}{h} + \frac{l_{0,2}-l_{0,3}}{h}\right) = \text{CFL} \in (0, 1]
\end{aligned}
\end{equation}
to get
\begin{equation}\label{eqn:A20}
 \alpha_0-\lambda_0  =    \frac{1}{2}(l_{0,2} + l_{0,3}) - \frac{l_{0,2} \, \lambda }{1+l_{0,2}-l_{0,3}}
\end{equation}
where we use $\lambda$ for the CFL number. Clearing denominators and simplifying, this is $\ge 0$ if and only if
\begin{equation}
    l_{02}^2-l_{03}^2 + l_{02} \, (1-2\lambda) + l_{03}  {\geq}0 .
\end{equation}
The left hand side is minimized when $\lambda = 1$, so we further simplify to checking if
\begin{align*}
    l_{02}^2-l_{03}^2 - l_{02}  + l_{03} & {\geq} 0 , 
\end{align*}
or equivalently, checking if
\begin{align}\label{eqn:A22}
    (l_{02}-l_{03})(l_{02}  + l_{03} - 1) & {\geq} 0 .
\end{align}
Equation \eqref{eqn:A22} is true, since we have $\alpha_0 = 1/2\,(l_{0,2}+l_{0,3})  \geq 1/2$ and $l_{0,2}\geq l_{0,3}$, thus showing the original equation \eqref{eqn:A20} is non-negative.
$\blacksquare$

\section{Stability of SRD for 1D Model Problem}\label{sec:gks}
This section proves a claim about stability from \cref{sec:simpleEx},with the help of a plot.

{\bf Theorem:} The upwind scheme \eqref{eq:pde} applied to the model problem \eqref{eqn:simpleMergeLeft} with one small cell with volume fraction $\alpha< 1/2$, merging left and without pre-merging,  is GKS stable for both original SRD and the new SRD weights.

We assume the interested reader is familiar with GKS theory \cite{GKS72}. 
We use the cell index notation and update formulas of \eqref{eqn:simpleMergeLeft}. The GKS theory for the discrete initial boundary value problem for dissipative schemes says to look for an eigenvalue $z$ with an eigenvector $\hat{U}$ that is bounded in space. If there are no solutions $|z|>1$ and no borderline cases with $|z|=1$ the scheme is GKS stable.

Thus we look for solutions of the difference equation of the form
$
     U^n_j = z^n \hat{U}_j .
$
The upwind scheme for the regular cells ($j \geq 1$ or $j \leq -2$) then satisfies
\begin{equation}\label{eqn:B1}
\begin{aligned}
    z \, \hat{U}_j &= \hat{U}_j (1-\lambda) + \lambda \hat{U}_{j-1} \\[.05in]
    \hat{U}_j &= \frac{\lambda}{(z-1+\lambda)} \, \hat{U}_{j-1} .
\end{aligned}
\end{equation}
This is a two term recurrence relation. In common GKS notation this may be written 
\begin{equation}\label{eqn:B2}
     \hat{U}_j = \left\{ 
     \begin{array} {ll} C_+ \cdot\kappa^j 
                   & \text{ if } j \geq \;1 \\[.3em]
                             C_- \cdot \kappa^j 
                   & \text{ if } j \leq -2 \;,
                             \end{array}
    \right. 
\end{equation}
with constants $C_-, C_+$ and 
\begin{equation}\label{eqn:B3}
         \kappa = \frac{\lambda}{z-1+\lambda} \;.
\end{equation}    
In the case $\left|\kappa\right|=1$ and $\kappa \neq 1$ then $\left|z\right|< 1$ 
because the scheme is von Neumann stable.
We analyze separately the remaining cases $\left|\kappa\right| < 1$ and $\left|\kappa\right|>1$.

In the case $\left|\kappa\right| > 1$ we show we must have $\left|z\right| < 1$.
Solving for $z$ instead of $\kappa$ in \eqref{eqn:B3} gives  
\begin{align*}
    z &= (1-\lambda) - \frac{\lambda}{\kappa} \\
    |z| &\leq (1-\lambda) + \frac{\lambda}{|\kappa|} .
\end{align*}
Since $|\kappa| > 1$ and $0 < \lambda \leq 1$ we immediately have $|z| < 1$.
Thus, there are no unstable ($\left|z\right|>1$) eigenvalues corresponding
to $\left|\kappa\right|>1$.

The case $\left|\kappa \right|<1$ is more complicated.
Observe that (B.2) implies that $C_-=0$ if $\hat{U}_j$ is bounded, since
$\left|C_-\kappa^j\right| \to \infty$ as $j \to - \infty$ . %
Thus, an eigenvector must have the form
\begin{align*}
     \hat{U}_j \; \;  &= 0 \;\text{ if } j\leq -2 \\
     \hat{U}_{-1} &= c \\
     \hat{U}_0  \; \; &= d  \\
     \hat{U}_j  \; \; &= C_+ \cdot \kappa^j \;\text{ if } j\geq 1 ,
\end{align*}
where we introduce constants $c$ and $d$ for the two cells where the recurrence relation doesn't apply.
We will examine the linear equations for $c$ and $d$ and show that they
have no unstable solutions.

\begin{figure}[h!]
    \centering
    \hfill
    \includegraphics[width=.4\linewidth,trim=10 0 30 10,clip]{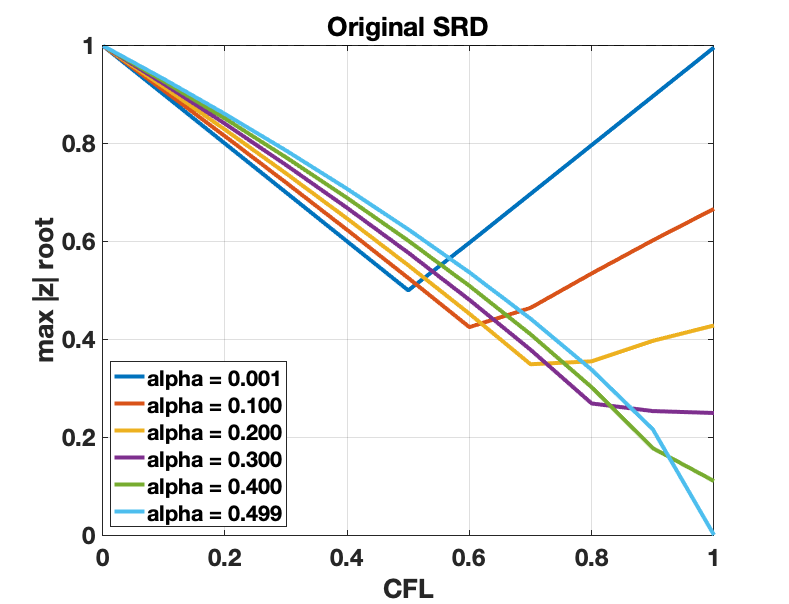}
    \hfill
    \includegraphics[width=.4\linewidth,trim=10 0 30 10,clip]{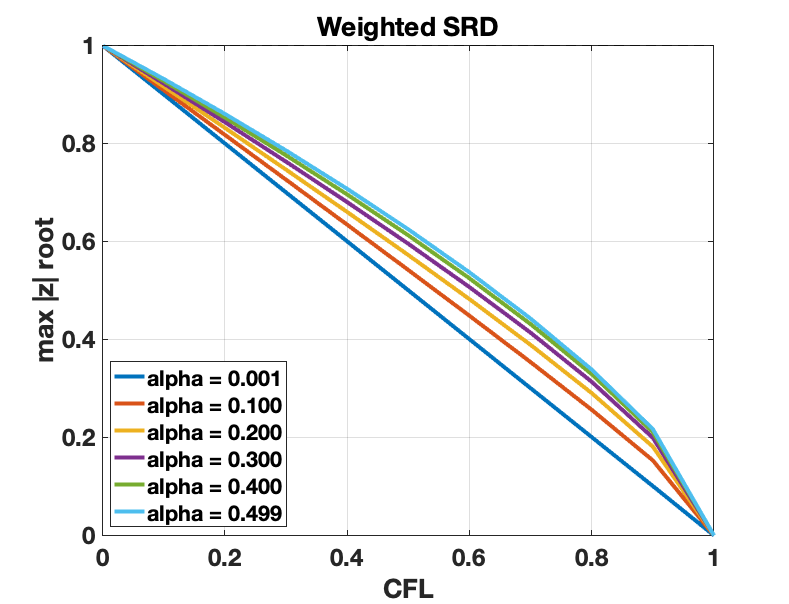}
    \hfill%
    \caption{Maximum root of GKS determinant is $\leq 1$ for the model problem with one small cell, merging left and without pre-merging. Left figure shows the original weights ($w=1/2$),  right shows the new weights ($w=\alpha$). All growth is bounded by 1, showing stability of  SRD although neither version is monotone.  \label{fig:gksfig}}
\vspace*{-.35in}
\end{figure}

The equations for $c$ and $d$ come from the update formulas for $U_0$ and $U_{-1}$ given in \eqref{eqn:simpleMergeLeft}. 
Writing $w$ instead of $w_{-1,-1}$, the formulas \eqref{eqn:simpleMergeLeft} become
\begin{equation}\label{eqn:gksmat}
\begin{aligned}
    z\,c &= \frac{(1 -w +\alpha \, w - \alpha\lambda w) }{\widehat V_0} \, c  + 
    \frac{(1-w)\,(\alpha-\lambda)}{ \widehat V_0} \, d \\[.07in]
    z\,d &=  \frac{(\lambda+(1-w)\,(1-\lambda))}{ \widehat V_0} \: c \: + \:  \frac{(\alpha-\lambda)}{\widehat V_0} \, d 
\end{aligned}
\end{equation}
In matrix form, this is
\begin{equation}\label{eqn:gksmat2}
\begin{pmatrix}
z - \frac{(1-w+\alpha w-\alpha\lambda w)}{\widehat V_0} & \frac{-(1-w)(\alpha-\lambda)}{\widehat V_0)}  \\[.08in]
\frac{-(\lambda+(1-w)(1-\lambda))}{\widehat V_0}  & z - \frac{(\alpha-\lambda)}{\widehat V_0} 
\end{pmatrix}
\begin{pmatrix}
    c \\[.08in]
    d
\end{pmatrix}
= 
\begin{pmatrix}
    0 \\[.08in]
    0
\end{pmatrix}
\end{equation}
The determinant of the matrix  is a quadratic polynomial in $z$.  The scheme is stable if there are no values of $|z |> 1$ where the determinant of \eqref{eqn:gksmat2} = 0.
Figure \cref{fig:gksfig} shows that all the roots have $\left|z\right|<1$ unless
$\lambda = 0$, which is the trivial case $\Delta t = 0$.
Note that this argument uses $\left|\kappa\right|<1$ only to 
establish that $C_-=0$ in (B.2).
Also, if $c=0$ and $d=0$ then also $C_+=0$ so $\hat{U}_j = 0$ for 
all $j$. 
This shows there are no non-zero eigenvectors corresponding to 
eigenvalues $z$ with $\left|z\right|\geq 1$ except the trivial 
case $z = 1$ and $\hat{U}_j \equiv 1$, which is GKS stability.

In the two cases of interest, the original SRD weights for the model problem set $w=1/2$. The new monotone weights set $w = \alpha$. \cref{fig:gksfig} plots the maximum absolute value of the two roots as a function of the CFL number $\lambda$, for several small cell fractions $\alpha < 0.5$. The values are clearly below 1, except the extremum value CFL = 1, which we check directly.  For the limiting case $\lambda = 1$, the root $z = \frac{2 \alpha-1}{2\alpha+1} < 1$  so for $\alpha > 0$ this is also stable.
The plots show the direction that $|z|$ moves in for the other values.

Summarizing, if $|\kappa| > 1$ there are no unstable $z, \kappa$ pairs at all. If $\kappa < 1$ there are unstable $z, \kappa$ pairs but they do not satisfy the boundary conditions. $\blacksquare$

\FloatBarrier

{\bf Acknowledgments}\label{sec:ack}
We thank the referees, whose suggestions have greatly improved this paper.
Mathematica \cite{Mathematica} was extremely valuable to double-check all the formulas in this paper.

\bibliographystyle{siamplain}
\bibliography{references}

\end{document}

%% file: ex_shared.tex
\usepackage{lipsum}
\usepackage{amsfonts}
\usepackage{graphicx}
\usepackage{epstopdf}
\usepackage{algorithmic}
\ifpdf
  \DeclareGraphicsExtensions{.eps,.pdf,.png,.jpg}
\else
  \DeclareGraphicsExtensions{.eps}
\fi

\newsiamremark{remark}{Remark}
\newsiamremark{hypothesis}{Hypothesis}
\crefname{hypothesis}{Hypothesis}{Hypotheses}
\newsiamthm{claim}{Claim}

\headers{Provably stable weighted SRD }{M. Berger, A. Giuliani}

\title{A new provably stable weighted state redistribution algorithm
}

\author{Marsha Berger\thanks{Flatiron Institute, New York City, New York 
  (\email{mberger@flatironinstitute.org}). Also Professor Emeritus, Courant Institute, New York University}
\and Andrew Giuliani\thanks{Flatiron Institute, New York City, New York 
  (\email{agiuliani@flatironinstitute.org}).}
}

\usepackage{amsopn}

%% file: one_small_right.tex
  \begin{tikzpicture}

  \definecolor{airforceblue}{rgb}{0.36, 0.54, 0.66}
   
   \def \a{0.6}
   \def \h{1}

   \draw (-3*\h , 0.0) -- (5*\h ,0.0);
   \draw[thick] (-1.5\h , 0.0) -- (3.\h,0.0);
   
   \draw[thick,<->, red!70!black] (0.5*\h, -0.8) -- (1.5*\h +\h*\a+0.05*\h,-0.8)  node[draw=none,fill=none,midway,below] {$\hat C_{0}$};

   \begin{scope}[yshift=.1cm]
   \draw[thick,<->, airforceblue] (-1.5*\h +0.05*\h, 0.8)  --  (-0.5*\h- 0.05*\h,0.8) node[draw=none,fill=none,midway,above] {$\hat C_{\sm 2}$};
   \draw[thick,<->, airforceblue] (-0.5*\h, 0.8)  --  ( 0.5*\h,0.8) node[draw=none,fill=none,midway,above] {$\hat C_{\sm 1}$};
   \draw[thick,<->, airforceblue] (1.5*\h +\h*\a-0.05*\h, 0.8)  --  (0.5*\h +\h*\a,0.8) node[draw=none,fill=none,midway,above] {$\hat C_{1}$};
   \draw[thick,<->, airforceblue] (2.5*\h +\h*\a, 0.8)  --  (1.5*\h +\h*\a+0.05*\h,0.8) node[draw=none,fill=none,midway,above] {$\hat C_{2}$};
   \end{scope}
   
   \draw[thick] (-1.5*\h ,-0.1) -- (-1.5*\h ,0.1);
   \draw[thick] (-0.5*\h,-0.1) -- (-0.5*\h,0.1);
   \draw[thick] ( 0.5*\h,-0.1) -- ( 0.5*\h,0.1);
   \draw[thick] (0.5*\h +\h*\a,-0.1) -- (0.5*\h +\h*\a,0.1);
   \draw[thick] (1.5*\h +\h*\a,-0.1) -- (1.5*\h +\h*\a,0.1);
   \draw[thick] (2.5*\h +\h*\a,-0.1) -- (2.5*\h +\h*\a,0.1);
   
  \node at (0.15*\h+\h*\a,0.4){\scriptsize $w_{0,0}$};
  \node at (\h + \h*\a,0.6){\scriptsize $w_{1,1}$};
  \node at (\h + \h*\a,0.2){\scriptsize $w_{1,0}$};
   
  \node at (-1*\h,-0.4){$h$};
   \node at (0, -0.4){$h$};
   \node at (0.15*\h+\h*\a,-0.4){$\alpha h$};
   \node at (\h + \h*\a,-0.4){$h$};
   \node at (2*\h + \h*\a,-0.4){$h$};
   
    \node at (4.25*\h ,-0.4){cell sizes};

  \end{tikzpicture}

%% file: 2D_cut_square/2D_1_main.tex
  \begin{tikzpicture}
  \tikzset{
  reverseclip/.style={insert path={(-4.2,-2.7) rectangle (-0.8,1.1)}}
}
\definecolor{beaublue}{rgb}{0.74, 0.83, 0.9};

\begin{scope}
\clip [draw](-4.,-1.75)--(-1,-1.1) --(-1,1) -- (-4.,1)-- (-4., -1.75); 

\draw (-5,-5) grid (5,5);
\draw [ultra thick, draw=green!70!black, fill=green!70!black, fill opacity=0.2] (-3,-2) rectangle  (-2,0);
\node (A) at (-2.5,-0.5) {$C_{1}$};

\node (A) at (-3.5,-0.5) {$C_{5}$};

\node (A) at (-1.5,-0.5) {$C_{4}$};
\end{scope}

\begin{scope}
\clip [reverseclip](-4.1,-1.77)--(-0.9,-1.08) --(-1.1,1) -- (-4.1,1)-- (-4.1, -1.75); 

\draw [ultra thick, draw=beaublue, fill=beaublue, fill opacity=0.2] (-2,-2) rectangle node[text=black, opacity=1]{}  (-1,-1);
\draw [ultra thick, draw=beaublue, fill=beaublue, fill opacity=0.2] (-3,-2) rectangle node[text=black, opacity=1]{}  (-2,-1);
\draw [ultra thick, draw=beaublue, fill=beaublue, fill opacity=0.2] (-4,-2) rectangle node[text=black, opacity=1]{}  (-3,-1);

\draw[<-]  (-3.5, -1.75)--  (-3.5,-2.2)  node [below] { $C_2$ };
\node (A) at (-2.5,-1.75) {$C_{0}$};
\node (A) at (-1.5,-1.6) {$C_{3}$};

\draw [ultra thick, draw=beaublue, fill=beaublue, fill opacity=0.2] (-3,-1) rectangle node[text=black, opacity=1]{ }  (-2,0);

\draw[->, thick]  (-2.5, -2.35)--  (-1.5,-2.12)  node [midway, below] { $\mathbf a$ };

\end{scope}
\node (A) at (-2.5,1.4) {Configuration 1};

  \end{tikzpicture}

%% file: 2D_cut_square/2D_2_main.tex
  \begin{tikzpicture}
  \tikzset{
  reverseclip/.style={insert path={(-4.2,-2.7) rectangle (-0.8,1.1)}}
}
\definecolor{beaublue}{rgb}{0.74, 0.83, 0.9};
\def\a{1};
\def\b{-1.25};

\begin{scope}
\clip [draw](-4.,-1.75 + \a + \b)--(-1,-1.1 + \a) --(-1,1) -- (-4.,1)-- (-4., -1.75 + \a+ \b); 

\draw (-5,-5) grid (5,5);
\draw [ultra thick, draw=green!70!black, fill=green!70!black, fill opacity=0.2] (-3,-2) rectangle  (-2,0);
\node (A) at (-2.5,-0.5) {$C_{1}$};

\node (A) at (-3.5,-0.5) {$C_{5}$};
\end{scope}

\begin{scope}
\clip [reverseclip](-4.1,-1.81 + \a+ \b)--(-0.8,-0.97+ \a) --(-1,1) -- (-4.1,1)-- (-4., -1.81 + \a); 

\node[] at (-1.45,-0.72) { $C_4$};

\draw [ultra thick, draw=beaublue, fill=beaublue, fill opacity=0.2] (-3,-1) rectangle node[text=black, opacity=1]{ }  (-2,0);

\draw [ultra thick, draw=beaublue, fill=beaublue, fill opacity=0.2] (-2,-2) rectangle node[text=black, opacity=1]{}  (-1,-1);
\draw [ultra thick, draw=beaublue, fill=beaublue, fill opacity=0.2] (-2,-1) rectangle node[text=black, opacity=1]{}  (-1,0);
\draw [ultra thick, draw=beaublue, fill=beaublue, fill opacity=0.2] (-3,-2) rectangle node[text=black, opacity=1]{}  (-2,-1);
\draw [ultra thick, draw=beaublue, fill=beaublue, fill opacity=0.2] (-4,-2) rectangle node[text=black, opacity=1]{}  (-3,-1);

\draw [ultra thick, draw=beaublue, fill=beaublue, fill opacity=0.2] (-3,-1) rectangle node[text=black, opacity=1]{ }  (-2,0);

\draw[<-]  (-3.5, -1.75)--  (-3.5,-2.2)  node [below] { $C_2$ };
\node (A) at (-2.5,-1.5) {$C_{0}$};
\node (A) at (-1.5,-1.5) {$C_{3}$};
\draw [ultra thick, draw=beaublue, fill=beaublue, fill opacity=0.2] (-3,-1) rectangle node[text=black, opacity=1]{ }  (-2,0);

\draw[->, thick]  (-2.5, -2.35)--  (-1.6,-1.7)  node [midway, below] { $\mathbf a$ };

\end{scope}
\node (A) at (-2.5,1.4) {Configuration 2};

  \end{tikzpicture}

%% file: 2D_cut_square/2D_7_main.tex
  \begin{tikzpicture}
  \tikzset{
  reverseclip/.style={insert path={(-4.2,-2.) rectangle (-0.8,1.1)}}
}
\definecolor{beaublue}{rgb}{0.74, 0.83, 0.9};
\def\a{1.75};
\def\b{-2.25};

\begin{scope}
\clip [draw](-3.75, -2.)--(0,1.65) --(0,2)  -- (-4.,2)-- (-4., -2.) -- (-3.75, -2.); 

\draw (-5,-2) grid (5,5);
\draw [ultra thick, draw=green!70!black, fill=green!70!black, fill opacity=0.2] (-3,-1) rectangle  (-1,1);
\node (A) at (-2.6,-0.5) {$C_{2}$};
\node (A) at (-2.5, 0.5) {$C_{5}$};
\node (A) at (-1.6, 0.5) {$C_{1}$};

\end{scope}

\begin{scope}
\clip [reverseclip](-3.75, -2.)--(0, 1.65)--(0, 2)   -- (-4.,2) ;
\node[] at (-1.45,-0.5) { $C_0$};

\draw [ultra thick, draw=beaublue, fill=beaublue, fill opacity=0.2] (-2,-1) rectangle node[text=black, opacity=1]{}  (-1,0);
\draw [ultra thick, draw=beaublue, fill=beaublue, fill opacity=0.2] (-3,-1) rectangle node[text=black, opacity=1]{}  (-2,0);
\draw [ultra thick, draw=beaublue, fill=beaublue, fill opacity=0.2] (-2,0) rectangle node[text=black, opacity=1]{}  (-1,1);

\end{scope}

\draw[->, thick]  (-1.85-.6, -1.85+.2)--  (-1-.6,-1+.2)  node [midway, below] { $\mathbf a$ };

\node (A) at (-2.1,2.3) {2-by-2 central merging};
  \end{tikzpicture}

%% file: 2D_cut_square/2D_6_main.tex
  \begin{tikzpicture}
  \tikzset{
  reverseclip/.style={insert path={(-4.2,-2.) rectangle (0.1,1.1)}}
}
\definecolor{beaublue}{rgb}{0.74, 0.83, 0.9};
\def\a{1.75};
\def\b{-2.25};

\begin{scope}
\clip [draw](-3.75, -2.)--(0,1.65) --(0,2)  -- (-4.,2)-- (-4., -2.) -- (-3.75, -2.); 

\draw (-5,-2) grid (5,5);
\draw [ultra thick, draw=green!70!black, fill=green!70!black, fill opacity=0.2] (-3,-2) rectangle  (0,1);
\node (A) at (-2.6,-0.5) {$C_{2}$};
\node (A) at (-2.5, 0.5) {$C_{5}$};
\node (A) at (-1.6, 0.5) {$C_{1}$};

\end{scope}

\begin{scope}
\clip [reverseclip](-3.75, -2.)--(0, 1.65)--(0, 2)   -- (-4.,2) ;

\draw [ultra thick, draw=beaublue, fill=beaublue, fill opacity=0.2] (-3,-2) rectangle node[text=black, opacity=1]{}  (-2,-1);

\draw [ultra thick, draw=beaublue, fill=beaublue, fill opacity=0.2] (-2,-1) rectangle node[text=black, opacity=1]{}  (-1,0);

\draw [ultra thick, draw=beaublue, fill=beaublue, fill opacity=0.2] (-1,0) rectangle node[text=black, opacity=1]{}  (0,1);

\draw [ultra thick, draw=beaublue, fill=beaublue, fill opacity=0.2] (-2,0) rectangle node[text=black, opacity=1]{}  (-1,1);

\draw [ultra thick, draw=beaublue, fill=beaublue, fill opacity=0.2] (-3,-1) rectangle node[text=black, opacity=1]{}  (-2,0);

\end{scope}
\node[] at (-0.45, 0.5) { $C_4$};
\node[] at (-1.45,-0.5) { $C_0$};
\node[] at (-2.45,-1.5) { $C_6$};

\draw[->, thick]  (-1.85-.3, -1.85+.2)--  (-1-.3,-1+.2)  node [midway, below] { $\mathbf a$ };

\node (A) at (-2.1,2.3) {3-by-3 central merging};

  \end{tikzpicture}

%% file: 2D_cut_square/2D_1.tex
  \begin{tikzpicture}
  \tikzset{
  reverseclip/.style={insert path={(-4.2,-2.7) rectangle (-0.8,1.1)}}
}
\definecolor{beaublue}{rgb}{0.74, 0.83, 0.9};

\begin{scope}
\clip [draw](-4.,-1.75)--(-1,-1.1) --(-1,1) -- (-4.,1)-- (-4., -1.75); 

\draw (-5,-5) grid (5,5);
\draw [ultra thick, draw=green!70!black, fill=green!70!black, fill opacity=0.2] (-3,-2) rectangle  (-2,0);
\node (A) at (-2.5,-0.5) {$C_{1}$};
\node (A) at (-3.5,-0.5) {$C_{5}$};
\node[] at (-1.5,-0.5) { $C_4$};
\end{scope}

\draw[<->] (-3.1, -1.55) -- (-3.1, -1) node [midway, left] { \footnotesize $l_{0, 2}$};
\draw[<->] (-2.1, -1.34) -- (-2.1, -1)  node [midway, left] {\footnotesize $l_{0, 3}$};

\begin{scope}
\clip [reverseclip](-4.1,-1.77)--(-0.9,-1.08) --(-1.1,1) -- (-4.1,1)-- (-4.1, -1.75); 

\draw [ultra thick, draw=beaublue, fill=beaublue, fill opacity=0.2] (-2,-2) rectangle node[text=black, opacity=1]{}  (-1,-1);
\draw [ultra thick, draw=beaublue, fill=beaublue, fill opacity=0.2] (-3,-2) rectangle node[text=black, opacity=1]{}  (-2,-1);
\draw [ultra thick, draw=beaublue, fill=beaublue, fill opacity=0.2] (-4,-2) rectangle node[text=black, opacity=1]{}  (-3,-1);

\draw[<-]  (-3.5, -1.75)--  (-3.5,-2.2)  node [below] { $C_2$ };
\node (A) at (-2.5,-1.75) {$C_{0}$};
\node (A) at (-1.5,-1.6) {$C_{3}$};

\draw [ultra thick, draw=beaublue, fill=beaublue, fill opacity=0.2] (-3,-1) rectangle node[text=black, opacity=1]{ }  (-2,0);

\draw[->, thick]  (-2.5, -2.35)--  (-1.5,-2.12)  node [midway, below] { $\mathbf a$ };

\end{scope}
\node (A) at (-2.5,1.4) {Configuration 1};

  \end{tikzpicture}

%% file: 2D_cut_square/2D_2.tex
  \begin{tikzpicture}
  \tikzset{
  reverseclip/.style={insert path={(-4.2,-2.7) rectangle (-0.8,1.1)}}
}
\definecolor{beaublue}{rgb}{0.74, 0.83, 0.9};
\def\a{1};
\def\b{-1.25};

\begin{scope}
\clip [draw](-4.,-1.75 + \a + \b)--(-1,-1.1 + \a) --(-1,1) -- (-4.,1)-- (-4., -1.75 + \a+ \b); 

\draw (-5,-5) grid (5,5);
\draw [ultra thick, draw=green!70!black, fill=green!70!black, fill opacity=0.2] (-3,-2) rectangle  (-2,0);
\node (A) at (-2.5,-0.5) {$C_{1}$};
\node (A) at (-3.5,-0.5) {$C_{5}$};

\end{scope}

\draw[<->] (-3.1, -1.4) -- (-3.1, -1) node [midway, left] { \footnotesize $l_{0, 2}$};
\draw[<->] (-1.9, -0.65) -- (-1.9, 0) ;
\node (A) at (-1.65,-0.25) {\footnotesize $l_{1,4}$};
\begin{scope}
\clip [reverseclip](-4.1,-1.81 + \a+ \b)--(-0.8,-0.97+ \a) --(-1,1) -- (-4.1,1)-- (-4., -1.81 + \a); 

\node[] at (-1.45,-0.72) { $C_4$};

\draw [ultra thick, draw=beaublue, fill=beaublue, fill opacity=0.2] (-3,-1) rectangle node[text=black, opacity=1]{ }  (-2,0);

\draw [ultra thick, draw=beaublue, fill=beaublue, fill opacity=0.2] (-2,-2) rectangle node[text=black, opacity=1]{}  (-1,-1);
\draw [ultra thick, draw=beaublue, fill=beaublue, fill opacity=0.2] (-2,-1) rectangle node[text=black, opacity=1]{}  (-1,0);
\draw [ultra thick, draw=beaublue, fill=beaublue, fill opacity=0.2] (-3,-2) rectangle node[text=black, opacity=1]{}  (-2,-1);
\draw [ultra thick, draw=beaublue, fill=beaublue, fill opacity=0.2] (-4,-2) rectangle node[text=black, opacity=1]{}  (-3,-1);

\draw [ultra thick, draw=beaublue, fill=beaublue, fill opacity=0.2] (-3,-1) rectangle node[text=black, opacity=1]{ }  (-2,0);

\draw[<-]  (-3.5, -1.75)--  (-3.5,-2.2)  node [below] { $C_2$ };
\node (A) at (-2.5,-1.5) {$C_{0}$};
\node (A) at (-1.5,-1.5) {$C_{3}$};
\draw [ultra thick, draw=beaublue, fill=beaublue, fill opacity=0.2] (-3,-1) rectangle node[text=black, opacity=1]{ }  (-2,0);

\draw[->, thick]  (-2.5, -2.35)--  (-1.6,-1.7)  node [midway, below] { $\mathbf a$ };
\end{scope}
\node (A) at (-2.5,1.4) {Configuration 2};

  \end{tikzpicture}